\documentclass[10pt]{article}
\usepackage{amstext}
\usepackage{amssymb}
\usepackage{amsmath}
\usepackage{latexsym}
\usepackage{euscript}
\usepackage[latin1]{inputenc}
\usepackage{epsfig}
\usepackage{enumerate}
\title{Minimal surfaces with helicoidal ends}
\author{Leonor Ferrer  \thanks{Research partially
supported by MCYT-FEDER grant number BFM2001-3489. \newline
 }
\and Francisco Martín $^{*}$ }
\date{}

\def\centerbmp#1#2#3{\vskip#2\relax\centerline{\hbox to#1{\special
  {bmp:#3 x=#1, y=#2}\hfil}}}
\def\r{\mathbb{R}}
\def\q{\mathbb{Q}}
\def\n{\mathbb{N}}
\def\c{\mathbb{C}}

\def\d{\mathbb{D}}

\newenvironment{proof}{\trivlist
\item[\hskip\labelsep{\em Proof}\,:]}{\hfill{$\Box$}\endtrivlist}

\def\ri{{\rm i}}
\def\re{{\rm e}}
\def\sh{\mathsf{h}}
\def\sd{\mathsf{d}}
\def\ts{\mathsf{t}}

\newcommand{\ee}{\begin{equation}}
\newcommand{\fe}{\end{equation}}

\newtheorem{lemma}{Lemma}
\newtheorem{remark}{Remark}
\newtheorem{theorem}{Theorem}

\begin{document}
\maketitle

%introduccion 
\section{Introduction and preliminaries} \label{sec:intro}
\begin{figure}[htpb]
	\begin{center}
		\includegraphics[width=.60\textwidth]{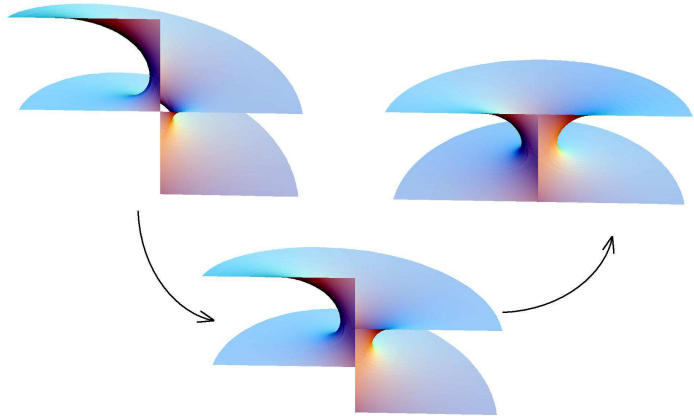}
	\end{center}
	\caption{The deformation that connects a fundamental piece of ${\cal H}_1$ to the López-Ritoré-Wei's surface.}
	\label{fig:deformation}
\end{figure}

In the last few years the study of minimal surfaces with helicoidal ends has gathered new speed. This is particularly the merit of D. Hoffman, H. Karcher and F. Wei who constructed the first examples of this kind of surfaces different from the helicoid. One of the examples constructed by these authors  was the so called singly-periodic genus-one helicoid, \cite{hkw}, that we will represent as ${\cal H}_1$. The helicoid ${\cal H}_1$ belongs to a continuous family of twisted periodic helicoids with handles that converges to a genus one helicoid. The continuity of this family of surfaces and the subsequent embeddedness of the genus one helicoid was obtained by D. Hoffman, M. Weber and M. Wolf in \cite{hww, weber}. In \cite{hww} the authors made a careful study of ${\cal H}_1$, with a new approach to the period problem associated to this surface. A thoughtful reading of this new approach establish a close relationship between ${\cal H}_1$ and an immersed minimal surface with planar ends, constructed by F.J. López, M. Ritoré and F. Wei in  \cite{l-r-w}. To be more precise, one observe that a fundamental piece of ${\cal H}_1$ can be obtained by deforming a fundamental piece of López-Ritoré-Wei's surface. The deformation consists of moving one of the connected components of the boundary of the surface following a vertical translation (see Fig. \ref{fig:deformation}). 
López and the second author \cite{l-m2, l-m1} constructed a family of minimal surfaces with planar ends based on López-Ritoré-Wei's, by modifying the angle between the horizontal boundary lines of the fundamental piece. 
\begin{figure}[h]
	\begin{center}
		\includegraphics[width=0.40\textwidth]{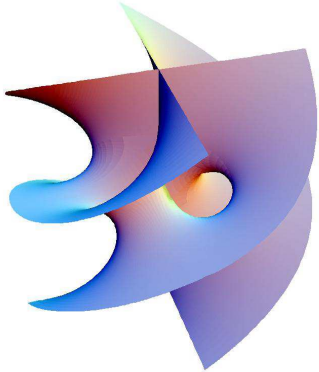}
	\end{center}
	\caption{A non-orientable example for angle $\pi/2$.}
	\label{fig:heli-2}
\end{figure}

Therefore, it is quite natural to ask whether is possible to construct new examples of minimal surfaces with helicoidal ends from the López-Martín examples. The main objective of the present paper is to describe this general deformation that connects López-Martín examples with a family of complete minimal surfaces with helicoidal ends that contains ${\cal H}_1$. 

These new surfaces, except for ${\cal H}_1$, are not embedded. However, if the angle between horizontal lines is $\frac{\pi}{n}$, with $n \in \n$, $n \geq 2$, then the only self-intersection of our surfaces occurs along the vertical axis. Furthermore, if $n$ is even the examples are non-orientable in $\r^3$. 

A simple proof of the embeddedness of the fundamental piece of our surfaces is obtained from the study of the above mentioned deformation. In the particular case of angle $\pi$ this argument provides another proof of the embeddedness for ${\cal H}_1$ (see Theorem \ref{embebi}). We also obtain the following uniqueness result
\begin{quote} \em Any complete, periodic, minimal surface containing a vertical line, whose quotient by vertical translations has genus one, contains two parallel horizontal lines, has two helicoidal ends and total curvature $-8 \pi$ is ${\cal H}_1$.
\end{quote}
This result was essentially obtain in \cite[Theorem 1]{hkw}. Our contribution consists of giving a new approach to the proof of the uniqueness of the period problem (see Remark \ref{hkw}).

The paper is lay out as follows. In Sect. \ref{deter} we determine the underlying complex structure and the Weierstrass data of a minimal disk bounded by a polygonal curve as in Fig. \ref{fig:gamma}. Thus we obtain a three-parameter family of Weierstrass data.  Sect. \ref{completas} is devoted to prove that this family of meromorphic data must contain, for each angle, an example with $q_1^-=q_2^+$ (see Fig. \ref{fig:gamma}). When the angle is a rational multiple of $\pi$, the surfaces obtained by successive Schwarz reflections about the straight lines are  complete and proper in $\r^3$. Finally, Sect. \ref{detalles} contains the technical details about the geometric functions that appear in Sect. \ref{completas}.
%datos de Weierstrass

\section{Determination of the Weierstrass representation}\label{deter}
As we announced in Sect.  \ref{sec:intro}, the fundamental piece of the minimal surfaces we wish to construct belongs to a family of minimal disks obtained moving one of the vertical segments upward in the López-Martín examples. Therefore, we are interested in the construction of properly embedded minimal disks whose boundary $\Gamma_{\beta \sd \sh}$ consists of the following configuration of straight lines:
\begin{figure}[ht]
	\begin{center}
		\includegraphics[width=.66\textwidth]{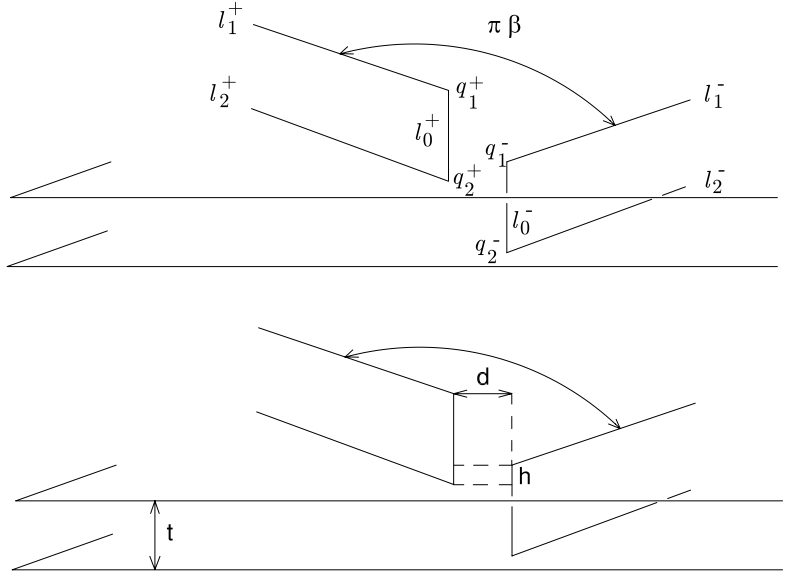}
	\end{center}
\caption{The curve $\Gamma_{\beta \sd \sh}$}\label{fig:gamma}
\end{figure}

Fix $\beta \in \: ]0,1]$, $\sd \geq 0$, $\sh \in \r$. Consider $\Pi$ a half-plane in $\r^3$ and denote by $\ell$ the boundary line of $\Pi$. Let $\ell_0$ be a line in $\Pi$ parallel to $\ell$ and $q_1^-$ and  $q_2^-$ two points in $\ell_0$. Denote by $\ell_0^-$ the segment $[q_1^-,q_2^-]$ and define $\ell_i^-$ as the half-line on $\Pi$ orthogonal to $\ell_0$ starting at $q_i^-, i=1,2$. Finally, we label $\ell_i^+$ and $q_j^+$ as the image of $\ell_i^-$ and $q_j^-$ for $i=0,1,2$, $j=1,2$, respectively, by a screw motion of axis $\ell$, angle $\beta \pi$ and vector $\ts \cdot \vec{n}$, $\ts \geq 0$, where $\vec{n}=\frac{\overrightarrow{q_2^- q_1^-}}{\big\|\overrightarrow{q_2^- q_1^-}\big\|}$. Write  $\sd=\mbox{dist}(\ell_0^+,\ell_0^-)$ and $\sh=\langle \overrightarrow{q_1^- q_2^+} ,\vec{n}\rangle$, where $\langle \cdot , \cdot \rangle$ denotes the usual inner product of $\r^3$.

Denote $\Pi^\pm$ as the plane that contains $\ell_0^{\pm} \cup \ell_1^{\pm} \cup \ell_2^{\pm}$. Finally, we define 
$$ \Gamma_{\beta \sd \sh}^+=\bigcup_{i=0}^2 \left(\ell_i^+ \right), \quad 
\Gamma_{\beta \sd \sh}^-=\bigcup_{i=0}^2 \left( \ell_i^-\right) , \quad
 \Gamma_{\beta \sd \sh}=\Gamma_{\beta \sd \sh}^+
\cup \Gamma_{\beta  \sd \sh}^- \; .$$
Observe that if $\ts=0$ we have exactly the family of examples given by López and Martín. Taking into account the geometric and topological properties of the surfaces we are starting at, we have to construct a properly immersed minimal surface $X=(X_1,X_2,X_3):M \longrightarrow \r^3$ with the following assumptions:
\begin{enumerate}[{\sc  {A}.1}] 
\item $M$ is homeomorphic to the closed unit disk $\overline{D}$ minus two boundary points $E_1$ and $E_2$. 
\item $X(\partial(M))= \Gamma_{\beta \sd \sh}$.
\item The surface has a symmetry respect to the line contained in the bisector plane of $\Pi^+$ and $\Pi^-$ that intersects orthogonally $\ell$ at the point $\frac{q_1^- + q_2^+}{2}$. 
\item If $\beta \in ]0,1[$ then $X(M)$ lies in the convex hull, ${\cal E}(\Gamma_{\beta \sd \sh})$, of $\Gamma_{\beta \sd \sh}$. If $\beta=1$, $X(M)$ lies in one of the two half slab determined by the plane $\Pi^+=\Pi^-$ and the planes orthogonal to $\ell$ containing $\ell_1^+$ and $\ell_2^-$.
\item If $\ell_0^+ \cap \ell_0^- = \emptyset$, $X$ is an embedding. In case $\ell_0^+ \cap \ell_0^- \neq  \emptyset$, the maps $X|_{M- \partial (M)^+}$ and $X|_{ M- \partial (M)^-}$ are injective, where  $\partial (M)^+$ and $\partial (M)^-$ are the two connected components of $\partial (M)$.
\end{enumerate}
From now on, we use a set of Cartesian coordinates, such that the half-plane $\Pi$ coincides with the half-plane
$\{(x_1,x_2,x_3) \in \r^3 \: \mid \: \sin\left(\frac{\pi \beta}{2}\right) x_1+ \cos\left(\frac{\pi \beta}{2}\right)(x_2+\frac{\sd}{2})=0, x_2 \leq 0\}$ and $q_2^-$ and $q_1^-$ are the points $(0,-\frac{\sd}{2},0)$ and $(0,-\frac{\sd}{2},\ts-\sh)$, respectively.

Assuming the above conditions and using the arguments presented in Sect. 2.1 of \cite{f-mp} with minor changes, it is easy to prove that the boundary has the following behavior:
\begin{lemma} \label{borde}
Up to relabelings, the  minimal immersion $X:M \longrightarrow \r^3$ satisfies $X(\partial(M)^+)=\Gamma_{\beta \sd \sh}^+$, $X(\partial(M)^-)=\Gamma_{\beta \sd \sh}^-$, $X ^{-1}(\ell_1^+) \cup X^{-1}(\ell_1^-)$ diverges to $E_1$ and 
$X ^{-1}(\ell_2^+) \cup X^{-1}(\ell_2^-)$ diverges to $E_2$.
\end{lemma}
Henceforth, $g$ and $\Phi_3$ will denote the Weierstrass data of a immersion $X:M \to \r^3$ satisfying the preceding five assumptions.
\subsection{The underlying complex structure of $M$} \label{M}

 The conformal type of $M$ can be easily determined using a global result on conformal structure of properly immersed minimal surfaces by P. Collin, R. Kusner, W.H. Meeks and H. Rosenberg
(see \cite{CKMR}). From Theorem  3.1 of \cite{CKMR} we obtain that $M$ is parabolic and hence, taking into account the topological type of $M$, $M$ is conformally equivalent to the closed unit disk $\overline{\d}$ minus two boundary 
points $E_1$ and $E_2$, where the biholomorphism extends piecewise analytically to the boundary.

Next, we prove that the Gauss map and Weierstrass data extend continuously to the ends.
\begin{lemma} \label{lem:phi3}
The 1-form $\Phi_3$ extends meromorphically to the ends. Even more, it has simple poles at the ends, with imaginary residues.
\end{lemma}
\begin{proof}
Let $(U_i,z)$ be a coordinate chart verifying that $U_i$ is biholomorphic to the upper half disk $\d^+=\{z \in \c \; : \; |z|<1 \leq 1+\mbox{Im}(z)\}$, $z(E_i)=0$, $z(\gamma_i^+ \cap U_i)=\d^+ \cap \r^+$, and $z(\gamma_i^-\cap U_i)=\d^+ \cap \r^-$, $i=1,2$. We know that $X_3$ is a bounded harmonic function on $U_i$ such that $X_3|_{\gamma_i^-}=C_i$, where $C_i$ is a constant, and $X_3|_{\gamma_i^+}=C_i+\ts$. So, the function $X_3+\frac{\ts}{\pi} \arg(z)$ can be continuously extended  to $E_i$, $i=1,2$. Then the function $X_3+\ri X_3^{\ast}- \frac{\ri \ts}{\pi} \log(z)$ is a holomorphic function on $U_i$ that extends to $E_i$, and so $\Phi_3-\frac{\ri\ts }{\pi} \frac{dz}{z}$ is a holomorphic 1-form on $U_i$. This concludes the proof. 
\end{proof}
\begin{lemma} \label{lem:gauss}
The Gauss map also extends and it is vertical at the ends.
\end{lemma}
\begin{proof}
For the case $0<\beta<1$, the arguments used in \cite[Theorem 3.12]{l-m2} also work in this setting. The case $\beta=1$ can be treated as in \cite[Proposition 3]{f-mp}.

Since $X(M)$ is contained between two horizontal parallel planes, the second assertion follows.
\end{proof}
The surface $X(M)$ can be extended by $180^{\circ}$ rotation about its boundary lines, to a complete surface (without boundary) in $\r^3$. 
Label $\widetilde{X}:\widetilde{M} \rightarrow \r^3$ the complete minimal immersion obtained in this way, where $\widetilde M$ is the corresponding Riemann surface without boundary, and let $(\widetilde g, \widetilde \Phi_3)$ denote its Weierstrass representation. Let $\cal S$ be the isometry of $\widetilde M$ induced by the symmetry described in assumption A.4.
We also denote ${\cal S}_j^\pm$, $j=0,1,2$, as the isometry of $\widetilde M$ induced by the $180^{\circ}$ rotation about the straight line containing $\ell_j^\pm$. Observe that: 
\begin{itemize}
\item $\tau_0={\cal S}_0^+ \circ {\cal S}_0^-$ is a horizontal translation whose translation vector, $\vec{v}_0$, is orthogonal to $\ell_0^-$ and $\ell_0^+$. Furthermore, the length of $\vec{v}_0$ is $2 \, \sd$;
\item $\tau_1={\cal S}_1^+ \circ {\cal S}_2^+$ is a vertical translation of translation vector $\vec{v}_1=(0,0,2 \, (\ts-\sh))$;
\item $\tau_2={\cal S}_2^+ \circ {\cal S}_2^-$ is a screw motion about the $x_3$-axis of angle $2 \, \beta \, \pi$ and translation vector $\vec{v}_2=(0,0,2 \, \ts).$
\end{itemize}

Let $\cal G$ be the subgroup of $\mbox{Iso}(\widetilde{M})$ generated by $\{\tau_0, \tau_1, \tau_2 \}$.  As $\cal G$ acts freely and properly discontinuously on $\widetilde M$, then the quotient ${\cal T}=\widetilde{M} / {\cal G}$ is a Riemann surface. Observe that $(d \widetilde{g})/ \widetilde{g}$ and $\widetilde{\Phi_3}$ can be induced in the quotient. We label $(dg)/g$ and $\Phi_3$ as the induced one-forms.

Taking into account Lemma \ref{borde} it is not hard to see that $\cal T$ has the topology of a torus minus four points. Moreover, this torus consists of four copies of $M$: $M \cup {\cal S}_0^+(M)  \cup {\cal S}_2^+(M) \cup ({\cal S}_2^+ \circ {\cal S}_0^+)(M)$ where  the boundary is identified according to the symmetries in $\cal G$. The compact torus is labeled as $\overline{\cal T}$. Note that $(dg)/g$ and $\Phi_3$ extends meromorphically to $\overline{\cal T}$.

Now, we need to determine the underlying complex structure of $\overline{\cal T}$. In order to do this, we will consider the symmetry ${\cal A}={\cal S}_2^+ \circ {\cal S}_0^+$. This symmetry is induced by a $180^{\circ}$ rotation in $\r^3$ about the orthogonal  line to $\Pi^+$ passing through $q_2^+$. $\cal A$ is a holomorphic involution that can be induced in the quotient ${\cal T}$. It can be also extended to $\overline{\cal T}$. We label $A$ as the induced involution in $\overline{\cal T}$. Note that $A$ exactly fixes four points $\{[q_i^+],[q_i^-]\}_{i=1,2}$ where $[p]$ denotes the class in $\overline{\cal T}$ of a point $p \in \widetilde{M}$. Using Riemann-Hurwitz formula, it is straightforward to check that $\overline{\cal T}/\langle A \rangle$ is conformally equivalent to the Riemann sphere $\overline{\c}$. If we label $u: \overline{\cal T} \longrightarrow \overline{\c}$ the canonical projection, then $u$ is an elliptic function on $\overline{\cal T}$. Furthermore, the branch points of $u$ coincide with  $\{[q_i^+],[q_i^-]\}_{i=1,2}$. Up to a Möbius transformation we can assume that $u([q_1^+])=0$, $u([q_2^-])=\infty$ and $u([q_1^-])=r \in\r^+$. Moreover, we label $s=u([q_2^+])$. Hence, $\overline{\cal T}$ is conformally equivalent to the algebraic elliptic curve $\left\{(u,v) \in \overline{\c}\times \overline{\c}\: \bigm| \: v^2=u(u- r) (u- s)  \right\}$. 

Let $S,S_0^+, S_2^+: \overline{\cal T} \longrightarrow \overline{\cal T}$ be the maps induced by $\cal S$, ${\cal S}_0^+$ and ${\cal S}_2^+$, respectively. Observe that $S_0^+$ is an antiholomorphic involution that fixes the branch points of $u$. This means that $u \circ S_0^+=\overline{u}$ and so $s \in \r$. On the other hand, $S$ is a holomorphic involution verifying $S([q_1^+])=[q_2^-]$ and $S([q_1^-])=[q_2^+]$. Furthermore, as the fixed points of $S$ are not at the boundary then $u \circ S=\frac{k}{u}$, where $k <0$. Up to the change $u \mapsto \frac{u}{\sqrt{-k}}$, we can assume that $u \circ S=-\frac{1}{u}$. In particular $s= -\frac{1}{r}$. Summarizing, 
\begin{enumerate}[(a)]
\item $\overline{\cal T} \equiv \{(u,v) \in \overline{\c} \times \overline{\c} \: \mid \: v^2= u (u-r) (r u+1) \}$ ,
\item $S(u,v)=(-\frac{1}{u},\frac{v}{u^2})$, $S_0^+(u,v)=(\overline{u}, \overline{v})$, and $S_2^+(u,v)=(\overline{u}, -\overline{v})$.
\end{enumerate}
Next, we will write our torus in a new way which is more suitable for our computations. Consider for $\rho \in ]0,\pi[$  the following torus
$$\overline{\cal N}=\left\{(z,w) \in \overline{\c} \times \overline{\c} \bigm| w^2=z^4+1-2\: z^2 \cos \rho \right\} \; .$$
It is not difficult to see that the map $(u,v)=B(z,w):\overline{\cal N} \longrightarrow \overline{\cal T}$, given by
$$u(z,w)=\frac{\left( \re^{-\ri \,\frac{\rho}{2} } + z \right) \,\left( -\re^{\ri \,\frac{\rho}{2} } + z \right) \,\cos (\frac{\rho}{2} ) + 
    \left( 1 - \sin (\frac{\rho}{2} ) \right)\, w }{\cos (\frac{\rho}{2} )\, w + \left( -\re^{-\ri \,\frac{\rho}{2} } - z \right) \,\left( -\re^{\ri \,\frac{\rho}{2}} + z \right) \,\left( 1 - \sin (\frac{\rho}{2} ) \right) }\; ,$$
$$v(z,w)=\frac{\frac{-\ri }{2}\,{\left( -\ri  + \re^{\ri \,\frac{\rho}{2} } \right) }^6\,\left( \re^{\ri \,\frac{\rho}{2} } - z \right) \,\left( 1 + \re^{\ri \,\frac{\rho}{2} }\,z \right) \,\left( 1 + z^2 \right) }{\re^{4\,\ri \,\frac{\rho}{2} }\, \left( -3 + \cos (\rho) + 4\,\sin (\frac{\rho}{2} ) \right) \,{\left( \cos (\frac{\rho}{2} )\,{w} + \left( -1 + \sin (\frac{\rho}{2} ) \right) \,\left( -1 + z^2 - 2\,\ri \,z\,\sin (\frac{\rho}{2} ) \right)  \right) }^2}\; ,$$
is a biholomorphism.

\begin{figure}[h]
	\begin{center}
		\includegraphics[width=0.33\textwidth]{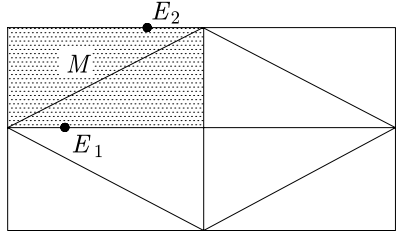}
	\end{center}
	\caption{The torus $\overline{\cal N}$ and the fundamental piece $M$.}
	\label{fig:recubridor}
\end{figure}

Note that the torus $\overline{\cal N}$ is a two-fold covering of the rhombic torus $\{(x,y) \in \overline{\c}^2 \; \; | \; y^2=x+ \frac{1}{x}-2 \cos \rho \}$. The covering map is given by $(z,w) \mapsto \left(z^2, \frac{w}{z} \right)$ (see Fig. \ref{fig:recubridor}). This family of rhombic tori (depending on $\rho$) coincides with those used by Hoffman, Karcher and Wei to construct the singly-periodic helicoids in \cite{hkw,hw}.

We will continue denoting by $S$, $S_0^+$ and $S_2^+$ the symmetries on the new torus $\overline{\cal N}$. According to (b) the expressions of these symmetries on $\overline{\cal N}$ are given by
\begin {equation} \label{newsim}S(z,w)=\left(\frac{1}{z},\frac{w}{z^2} \right)
 \; , \quad  S_0^+(z,w)=\left( \frac{1}{ \overline{z}},- \frac{\overline{w}}{ \overline{z}^2} \right) \; ,\quad  S_2^+(z,w)=(-\overline{z}, \overline{w})  \; .
\end{equation}
For the sake of brevity, when $z^4+1-2\: z^2 \cos \rho \in \r ^+$  we denote:
$$z_+=\left( z,+\sqrt{z^4+1-2\: z^2 \cos \rho}  \right)\; ,\quad z_-=\left(z,-\sqrt{z^4+1-2\: z^2 \cos \rho}  \right)\; .$$

Now we need to identify the punctures in this torus. Note that $S_2^+$ fixes the ends and $S_0^+$ interchanges them.
Taking  (\ref{newsim}) into account we deduce that the ends are $E=\{\ri a_+,-\ri a_-,\frac{\ri}{a}_+,-\frac{\ri}{a}_- \}$, where $a \in \r$. Up to relabeling, we can assume $a \in [0,1[$. We denote ${\cal N}=\overline{\cal N} - E$. 

On $\overline{\c}$ define the following set of curves:
$$s_0^+= \left\{{\rm e}^{{\rm i}\tfrac{t}{2}}\; \mid \; t \in \left[\rho,\pi \right] \right\} \; ,   
s_1^+=\{\lambda {\rm i} \; \mid \; \lambda \in \:]a,1]\} \; ,
s_2^+=\{\lambda {\rm i} \; \mid \; \lambda \in [1,\infty [ \,\cup \, ]-\infty,-\tfrac{1}{a}[\}\cup \{\infty \} \; ,$$
$$s_0^- = \left\{{\rm e}^{{\rm i}\tfrac{t}{2}}\; \mid \; t \in \left[-\pi,-\rho\right] \right\} \; , 
 s_1^-=\{\lambda {\rm i} \; \mid \; \lambda \in [-1,a [\} \; ,
 s_2^- =\{\lambda {\rm i} \; \mid \; \lambda \in \:]-\tfrac{1}{a},-1]\} \; .$$

Label $\gamma_i^+=z^{-1}(s_i^+)$, $\gamma_i^-=z^{-1}(s_i^-)$, $i=0,1,2$. In order to determine the domain in $\cal N$ that corresponds to our fundamental piece we observe that the set of fixed points of $S_0^+$ is $\gamma_0^+ \cup \gamma_0^- \cup S_2^+(\gamma_0^+ \cup \gamma_0^- )$ and the set of fixed points of $S_2^+$ is $\gamma_1^+ \cup \gamma_1^- \cup \gamma_2^+ \cup \gamma_2^-\cup S_0^+(\gamma_1^+ \cup \gamma_1^- \cup \gamma_2^+ \cup \gamma_2^-)$. According to this we can identify $M$ with the closure in ${\cal N}$ of the connected component of $z^{-1}(\overline\c - (\bigcup _{i=0}^2(s_i^+ \cup s_i^-)))$ containing the point $P_0=1_+$. We will label this domain as $M_{a,\rho}$ .
\begin{figure}[htbp]
	\begin{center}
		\includegraphics[width=0.33\textwidth]{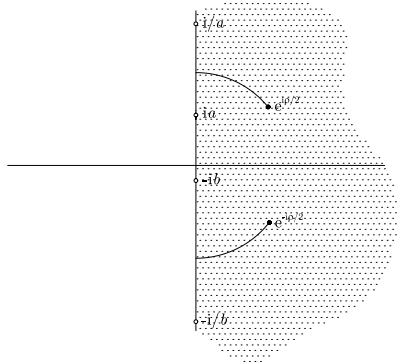}
	\end{center}
	\caption{The $z$-projection of the domain $M_{a,\rho}$.}
	\label{fig:darwin}
\end{figure}

Define $\gamma^+$ and $\gamma^-$ by:
$$ \gamma^+=\bigcup_{i=0}^2 \gamma_i^+, \quad \gamma^-=\bigcup_{i=0}^2 \gamma_i^-.$$ It is clear that $\partial (M)= \gamma^+ \cup \gamma^-$. Furthermore, 
note that $z|_{\gamma_i^+}$ and $z|_{\gamma_i^-}$ are bijective maps onto $s_i^+$ and $s_i^-$, respectively,
$i=1,2$. However, $\gamma_0^+$ and $\gamma_0^-$ consist of two copies of $s_0^+$ and $s_0^-$, respectively.

\subsection{The complex height differential}
According to Lemma \ref{lem:phi3} the height differential $\Phi_3$ has a simple pole with imaginary residue at the ends. As $\overline{\cal N}$ is a torus, then $\Phi_3$ has as many zeros as poles. Moreover, it is easy to see that the symmetries act on $\Phi_3$ as follows
\begin{equation} \label{simphi3}
S^*\Phi_3=-\Phi_3\;, \quad (S_0^+)^* \Phi_3= \overline{\Phi_3}\; , \quad (S_2^+)^* \Phi_3=- \overline{\Phi_3} \; .
\end{equation}
Both facts imply that $\Phi_3$ has four zeros of order one and they have this form $V=\{ - \ri b_+,\ri b_-,-\frac{\ri}{b}_+,\frac{\ri}{b}_- \}$, where $b \in ]-1,1[$. All this information lead us to:
\begin{equation} 
\Phi_3= \lambda \frac{w+c_2 \ri z}{w+c_1 \ri z} \frac{dz}{w} \; , 
\end{equation}
where $\lambda >0$ and
\begin{equation}\label{ces} c_1=\frac{1}{a} w({\rm i} a_+)= \frac{\sqrt{a^4+1+ 2\: a^2 \cos \rho}}{a}  \; , c_2=- \frac{1}{b} w(-{\rm i} b_+)= -\frac{\sqrt{b^4+1+ 2\: b^2 \cos \rho}}{b} \; . \end{equation}
\subsection{The Gauss map} \label{g}
The objective of this subsection is to find the expression of the Gauss map of our examples. From Lemma \ref{lem:gauss}, we have that the normal vector at the ends must be vertical. Then we can assume $g(\ri a_+)=0$. Furthermore, taking A.2 into account we deduce that the behavior of $g$ in a neighborhood of $\ri a_+$ is given by $g(z)= z^{\beta}$. Since $\Phi_3$ has zeros of order one at the points in $V$ we deduce that $g$ has at these points either simple poles or zeros of order one, in any case the points in $V$ are points where the normal vector is vertical. To obtain more information about the normal vector at the point  $- \ri b_+$ we need to return to our initial conditions. 

In fact, using {\sc A.4} and the interior maximum principle one can prove that $X(M)\cap \{x_3=0\}= \ell_2^-$, $X(M) \cap \{x_3=2 \ts-\sh\}= \ell_1^+ $, $X(M) \cap \{\sin \left( \frac{\pi \beta}{2} \right) \,x_1+ \cos \left( \frac{\pi \beta}{2} \right)\left(x_2+\frac \sd2 \right)=0\}= \Gamma_{\beta \sd \sh}^- $ and $X(M) \cap \{\sin \left( \frac{\pi \beta}{2} \right) \,x_1- \cos \left( \frac{\pi \beta}{2} \right)\left(x_2-\frac \sd2\right)=0\}= \Gamma_{\beta \sd \sh}^+ $. By studying the intersection of $X(M)$ with a horizontal plane containing $\ell_1^-$, we deduce the existence of a point with vertical normal at $\ell_1^-$. Clearly, this point must be $- \ri b_+ $. Suppose that $g$ has at the point  $- \ri b_+$ a simple pole. Then, as $X(- \ri b_+) \in \ell_1^-$ and $g(\ri a_+)=0$, should exists a point in $\ell_1^-$, different from the point $q_1^-$, whose tangent plane is $\{\sin \left( \frac{\pi \beta}{2} \right) \,x_1+ \cos \left( \frac{\pi \beta}{2} \right)\left(x_2+\frac \sd2 \right)=0\}$, but this contradicts what we have obtain previously. So, $g$ has at $- \ri b_+$ a zero of order one. 

 Finally, the behavior of the symmetries at the points in $E \cup V$, allows us to deduce that distribution of poles and zeros of the multivalued function $g$ using non-integral exponents must be as follows
\begin{center}
\begin{tabular}{|c|c|c|c|c|c|c|c|}
 $\ri a_+$ & $-\ri a_- $ & $\frac{\ri}{a}_+$ & $-\frac{\ri}{a}_-$ & $- \ri b_+$  & $\ri b_-$ & $-\frac{\ri}{b}_+$
& $\frac{\ri}{b}_-$ \\ \hline
$ $ & $ $ & $ $ & $ $ & $ $  & $ $ & $ $ & $ $ \\ 
$0^{\beta}$& $ 0^{\beta}$   & $\infty^{\beta}$ & $\infty^{\beta}$ & $0^1$& $0^1$& $\infty^1$ & $\infty^1$ \\ \hline  
\end{tabular}
\end{center}
Thus, $\frac{dg}{g}$ (recall that this is a meromorphic function on $\overline{\cal N}$), have only simple poles at the points in $E \cup V$ and the residues of  $\frac{dg}{g}$ at these points are given in the following table 
\begin{center}
\begin{tabular}{c|c c c c c c c c}
 $p$ & $\ri a_+$ & $-\ri a_- $ & $\frac{\ri}{a}_+$ & $-\frac{\ri}{a}_-$ & $- \ri b_+$  & $\ri b_-$ & $-\frac{\ri}{b}_+$
& $\frac{\ri}{b}_-$ \\ \hline

${\rm Residue}(\frac{dg}{g},p)$ & ${\beta}$& $ {\beta}$   & $-{\beta}$ & $-{\beta}$ & $1$& $1$& $-1$ & $-1$  
\end{tabular}
\end{center}
Furthermore, from the definitions of the symmetries in (\ref{newsim}) we have
$$ S^*(\Phi_1-\ri \Phi_2)=\re^{-\ri \beta \pi}(\Phi_1+\ri \Phi_2)\;, (S_0^+)^* (\Phi_1-\ri \Phi_2)= -(\overline{\Phi_1}-\ri \overline{\Phi_2})\; , (S_2^+)^*(\Phi_1-\ri \Phi_2)=\re^{-2\ri \beta \pi}(\overline{\Phi_1}+\ri \overline{\Phi_2}) \; .   $$
Now, taking into account (\ref{simphi3}) and that $g=\frac{\Phi_3}{\Phi_1-\ri \Phi_2}$ we infer
\begin{equation}\label{acciong}
S^{\ast}\left( \frac{dg}{g} \right)= -\frac{dg}{g} \; ,\quad (S_0^+)^{\ast} \left(\frac{dg}{g}\right)= -\overline{\left(\frac{dg}{g}\right)} \; , \quad (S_2^+)^{\ast} \left(\frac{dg}{g}\right)= \overline{\left(\frac{dg}{g}\right)} \; .
\end{equation}
All the facts above presented imply that $\frac{dg}{g}$ can be written as
$$\frac{dg}{g}= \eta_1+ a_3 \eta_2 \;,$$
where $a_3 \in \r$ and
$$\eta_1=\frac{{\rm i} \beta}{a_1} \frac{dz}{w+c_1{\rm i}z}+\frac{{\rm i}}{a_2} \frac{dz}{w+c_2{\rm i}z} \; , \quad 
\eta_2= {\rm i} \frac{dz}{w} \; ,$$ 
\begin{eqnarray} \label{a1} a_1 & = &{\rm Residue}\left(\frac{dz}{w+c_1{\rm i}z},{\rm i} a_+ \right)  = \frac{a \sqrt{a^4+1+2 a^2 \cos(\rho)}}{1-a^4} \; , \\ a_2 &=& {\rm Residue}\left(\frac{dz}{w+c_2{\rm i}z},-{\rm i} b_+ \right)  =  -\frac{b \sqrt{b^4+1+2 b^2 \cos(\rho)}}{1-b^4} \; . \label{a2} \end{eqnarray}

Now, we must prove that there exist $b(a,\rho,\beta)$ and $a_3(a,\rho,\beta)$ so that the Gauss map $g=\exp \left(\int _{1_+} \frac{dg}{g}\right)$ verifies on $M$ the other required conditions. Taking into account the symmetry $S$ it is sufficient to pay attention to $\{(z,w) \in M \mid {\rm Re}(z) \geq 0\}$. To translate these conditions into equations we need some terminology.

The curve $\gamma_0^+$ consists of two copies, $\delta_1$ and $\delta_2$,  of $s_0^+$. We can assume that $\delta_1(t)$ and $\delta_2(t)$ are the two lifts to $M$ of the curve $\re^{\ri \frac{t}{2}}$, $t \in [\rho,\pi]$, in the $z$-plane, satisfying $\delta_1(\pi) \in \gamma_1^+$ and $\delta_2(\pi) \in \gamma_2^+$, respectively. Let $\delta(t)$ be the lift to $M$  of the curve $\re^{\ri \frac{t}{2}}$, $t \in [0,\rho]$, in the $z$-plane. Observe $\delta(0)=1_+$ and $\delta(\rho)=\delta_1(\rho)=\delta_2(\rho).$ With this notation we have:
\begin{itemize}
\item As we wish that $|g|=1$ on $\gamma_0^+$ we have to impose
\begin{equation}{\rm Re}\: \left(\int_{\delta} \frac{dg}{g} \right)=0 \; . \label{eq1} \end{equation}
\item We also have to impose that $g(\ri _+)=g(\ri _-)$. Then we have the condition
\begin{equation} {\rm Im}\: \left(\int_{\widetilde{\delta}}\frac{dg}{g} \right)=0 \; ,  \label{eq2} \end{equation} 
where $\widetilde{\delta}=-\delta_1+ \delta_2$, with $-\delta_1(t)=\delta_1(\pi-t)$.
\end{itemize}
On ${\cal N}$ we consider the curves $\gamma_1$ and $\gamma_2$ below described:
\begin{itemize}
\item $\gamma_1$ is the curve in ${\cal N}$ given by $\gamma_1=\alpha_1 -(S_0^+)_{\ast}(\alpha_1)$, where $\alpha_1=-S_{\ast}(\delta)+ \delta$.
\item $\gamma_2$ is the curve in ${\cal N}$ given by $\gamma_2=-\alpha_2 +\alpha_3$, where $\alpha_2=\delta_1-(S_2^+)_{\ast}(\delta_1)$ and $\alpha_3=\delta_2-(S_2^+)_{\ast}(\delta_2)$.
\end{itemize}
Observe that $\{\gamma_1,\gamma_2 \}$ is a canonical homology base of $\overline{\cal N}$ and $\gamma_2=\widetilde{\delta} -(S_2^+)_{\ast}(\widetilde{\delta})$. Thus, from (\ref{acciong}) we have
$$\int_{\gamma_1} \frac{dg}{g}= \int_{\alpha_1} \frac{dg}{g} +\overline{ \int_{\alpha_1} \frac{dg}{g} }=2\:{\rm Re}\left(\int_{-S_{\ast}(\delta)+ \delta} \frac{dg}{g} \right)=4\:{\rm Re}\left(\int_{\delta} \frac{dg}{g} \right)\; , $$
$$ \int_{\gamma_2} \frac{dg}{g}= \int_{\widetilde{\delta}} \frac{dg}{g} -\overline{ \int_{\widetilde{\delta}} \frac{dg}{g}}=2\:{\rm Im}\left(\int_{\widetilde{\delta}} \frac{dg}{g} \right)\; . $$
Therefore the equations (\ref{eq1}) and (\ref{eq2}) are equivalent to the system:
\begin{equation} \int_{\gamma_i} \frac{dg}{g}= \int_{\gamma_i} \eta_1+a_3  \int_{\gamma_i} \eta_2=0 \; \quad {\rm for} \quad i=1,2 \; . \label{eq3}
\end{equation}
In order to solve (\ref{eq3}) it is sufficient to prove that there exists $b(a,\rho,\beta)$ such that
\begin{equation}\label{eq4} \det \left( \begin{array}{cc} \int_{\gamma_1} \eta_1 & \int_{\gamma_1} \eta_2 \\ \int_{\gamma_2}\eta_1 & \int_{\gamma_2}\eta_2  \end{array} \right)= \int_{\gamma_1} \eta_1 \int_{\gamma_2}\eta_2  -  \int_{\gamma_1} \eta_2 \int_{\gamma_2}\eta_1 =0 \, .
\end{equation}
Applying the bilinear relations of Riemann to the 1-forms $\eta_1$ and $\eta_2$ we obtain:
\begin{equation} \int_{\gamma_1} \eta_1 \int_{\gamma_2}\eta_2  -  \int_{\gamma_1} \eta_2 \int_{\gamma_2}\eta_1 = 2 \pi \ri \sum_{p \in E \cup V}
f(p) \:{\rm Residue}(\eta_1,p) \; , \label{eq5} \end{equation}
where $f$ is a primitive of $\eta_2$ on the simply connected domain of ${\cal N}$, $\Omega$, obtained by removing the curves $\gamma_1$ and $\gamma_2$ (see \cite{farkas}). We choose $f$ so that $f(\re^{\ri \frac{\rho}{2}})=0$. From (\ref{newsim}) we obtain that $(S_0^+)^{\ast}(\eta_2)=-\overline{\eta_2}$. Then we get
$$f(S_0^+(p))=\int_{\re^{\ri \frac{\rho}{2}}}^{S_0^+(p)}\eta_2=\int_{\re^{\ri \frac{\rho}{2}}}^{p}(S_0^+)^{\ast}(\eta_2)=-\overline{\int_{\re^{\ri \frac{\rho}{2}}}^{p}\eta_2}=-\overline{f(p)} \; .$$
On the other hand we can consider the holomorphic transformation on ${\cal N}$ given by $T(z,w)=(z,-w)$. As $T^\ast(\eta_2)=-\eta_2$ we also obtain
$$f(T(p))=\int_{\re^{\ri \frac{\rho}{2}}}^{T(p)}\eta_2=\int_{\re^{\ri \frac{\rho}{2}}}^{p}T^{\ast}(\eta_2)=-\int_{\re^{\ri \frac{\rho}{2}}}^{p}\eta_2=-f(p) \; .$$
Notice that in the above computations we have used that $S_0^+(\Omega)=T(\Omega)=\Omega.$

Therefore, the equality in (\ref{eq5}) can be written as
\begin{eqnarray} 
 \hspace{-0.3cm}\int_{\gamma_1} \eta_1 \int_{\gamma_2}\eta_2  -  \int_{\gamma_1} \eta_2 \int_{\gamma_2}\eta_1 &\hspace{-0.3cm}=& 
\hspace{-0.3cm}4 \:\pi \ri \:{\rm Re}\left(
\beta(f(\ri a_+)-f(-\ri a_+))+f(-\ri b_+)-f(\ri b_+)\right)= \nonumber\\
\hspace{-0.3cm}4 \pi \ri \:{\rm Re}\left(
\beta \int_{-\ri a_+}^{\ri a_+}\eta_2-\int_{-\ri b_+}^{\ri b_+}\eta_2\right)&\hspace{-0.3cm}=&\hspace{-0.3cm}-8 \pi  \ri \left( \beta \int_0^a\frac{dt}{\sqrt{t^4+1+2 t^2 \cos \rho}} -\int_0^b\frac{dt}{\sqrt{t^4+1+2 t^2 \cos \rho}}\right) \label{eq6}
\end{eqnarray}
Taking into account (\ref{eq6}), the equation (\ref{eq4}) is satisfied if $F(a,b,\rho,\beta)=0$, where
\begin{equation} F(a,b,\rho,\beta)=\beta \int_0^a\frac{dt}{\sqrt{t^4+1+2 t^2 \cos \rho}} -\int_0^b\frac{dt}{\sqrt{t^4+1+2 t^2 \cos \rho}}\; .\label{F} \end{equation}
The function $F$ can be expressed as
\begin{equation} F(a,b,\rho,\beta)= \int_0^1 \left( \frac{\beta a}{\sqrt{a^4 t^4+1+2 a^2 t^2 \cos \rho}} -\frac{b }{\sqrt{b^4 t^4+1+2 b^2 t^2 \cos \rho}}\right)dt\; .\label{F1} \end{equation}
Since $0 < \beta \leq 1$ we have that $F(a,a,\rho,\beta) \leq 0$. Moreover, it is easy to check that $\lim_{b \rightarrow 0} F(a,b,\rho,\beta) \geq 0$. Now, for $b \in [0,a]$ we have
$$\frac{\partial F}{\partial b}(a,b,\rho,\beta)=-\int_0^1\frac{(1-b^4 t^4) dt}{(b^4 t^4+1+2b^2 t^2 \cos \rho)^{\frac{3}{2}}} < 0 \; .$$
From the above settings, we infer that there exists a unique $b(a,\rho,\beta) \in \left[ 0,a \right]$ such that $F(a,b(a,\rho,\beta),\rho,\beta)=0$ and therefore verifies equation (\ref{eq4}). Moreover, one can give an explicit, but rather long, formula for the function $b$ in terms of elliptic functions and so $b$ is a real analytic function. Consequently, there exists $a_3(a,\rho,\beta) \in \r$ solution of the system (\ref{eq3}). From equation (\ref{eq1}) we obtain
\begin{equation} \label{a31} a_3(a,\rho,\beta)= - \frac{{\rm Re}\int_{\delta} \eta_1}{{\rm Re}\int_{\delta} \eta_2}=-\frac{\int_0^{\rho} \left( \frac{\beta}{a_1} \frac{1}{a^2+\frac{1}{a^2}+2 \cos(t)}+\frac{1}{a_2} \frac{1}{b^2+\frac{1}{b^2}+2 \cos(t)}\right) \sqrt{2(\cos(t)-\cos(\rho))}dt}{\int_0^{\rho} \frac{1}{\sqrt{2(\cos(t)-\cos(\rho))}}dt} \; .\end{equation}
Analogously, from (\ref{eq2}) we obtain an alternative definition for $a_3$ given by
\begin{equation} \label{a32} a_3(a,\rho,\beta)= -\frac{{\rm Im}\int_{\widetilde{\delta}} \eta_1}{{\rm Im}\int_{\widetilde{\delta}} \eta_2}=\frac{\int_{\rho}^\pi \left( \frac{\beta}{a_1} \frac{1}{a^2+\frac{1}{a^2}+2 \cos(t)}+\frac{1}{a_2} \frac{1}{b^2+\frac{1}{b^2}+2 \cos(t)}\right) \sqrt{2(\cos(\rho)-\cos(t))}dt}{\int_{\rho}^\pi \frac{1}{\sqrt{2(\cos(\rho)-\cos(t))}}dt} \; .\end{equation}
\section{The complete examples} \label{completas}
\begin{figure}[h]
	\begin{center}
		\includegraphics[width=0.33\textwidth]{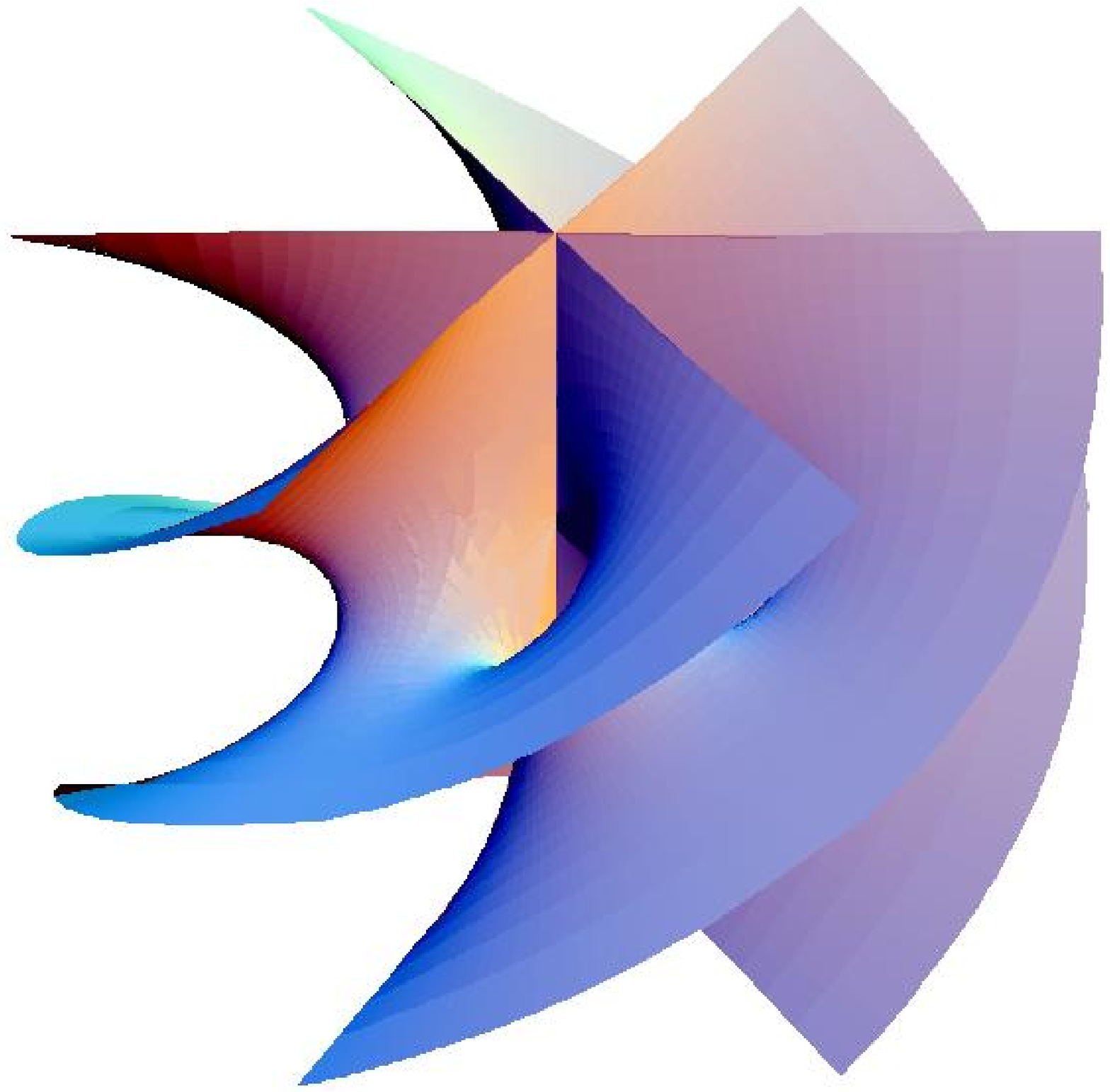}
	\end{center}
	\caption{A complete surface constructed with a fundamental piece of angle $\pi/3$.}
	\label{fig:heli-3}
\end{figure}

Our geometric assumptions of the minimal disk have led us to an explicit tree-parameter family  of Weierstrass data: For $a \in [0,1[$, $\rho \in \: ]0,\pi[$ and $\beta \in ]0,1]$, the disk $M_{a,\rho}$ defined at the end of paragraph \ref{M} and the meromorphic data
\begin{eqnarray} \label{Phi3} g=\exp \left( \int_{1_+}\frac{{\rm i} \beta}{a_1} \frac{dz}{w+c_1{\rm i}z}+\frac{{\rm i}}{a_2} \frac{dz}{w+c_2{\rm i}z}+{\rm i} a_3   \frac{dz}{w} \right) \; ,\;
\Phi_3= \lambda \frac{w+c_2 \ri z}{w+c_1 \ri z} \frac{dz}{w} \; ,   \end{eqnarray}
where $b$ is the function satisfying $F(a,b(a,\rho,\beta),\rho,\beta)=0$ with $F$ defined in (\ref{F}), $a_3$ is given by either (\ref{a31}) or (\ref{a32}), $\lambda >0$ and $c_i, a_i$ are given in (\ref{ces}), (\ref{a1}) and (\ref{a2}).

In addition, we can prove
\begin{theorem} \label{teo:grande}
A minimal immersion satisfying assumptions {\sc A.1}-{\sc A.5} has Weierstrass data of the form (\ref{Phi3}) with $a \in [0,1[$, $\rho \in \:]0,\pi[$ and $\beta \in ]0,1]$.

Conversely, for $a \in [0,1[$, $\rho \in \: ]0,\pi[$ and $\beta \in ]0,1]$ the Weierstrass data (\ref{Phi3}) define  a proper minimal immersion $X_{a,\rho,\beta}:M_{a,\rho} \longrightarrow \r^3$ that fulfills the assumptions {\sc A.1}-{\sc A.4}. Furthermore, 
$X_{|\gamma^+}$ and $X_{|\gamma^-}$ are injective.
\end{theorem}
\begin{proof}
The first part of the theorem is a direct consequence of the development along the previous section. For the sake of brevity, throughout the proof of the converse part of the theorem we denote $M=M_{a,\rho}$ and $X=X_{a,\rho,\beta}$. Clearly, a surface represented by the data as given in (\ref{Phi3}) on $M$ satisfies {\sc A.1} and {\sc A.3}. To see {\sc A.2} we parametrize the curves $\gamma_1^{\pm}$ as $\gamma_1^+(t)= \ri t_+ $, $t \in ]a,1]$  and $\gamma_1^-(t)= \ri t_+ $, $t \in [-1,a[$. At this point, we are interested in compute $g(\pm \ri_+)$. In order to do this, we consider the curve $\gamma=\alpha_1+\alpha_2-(S_2^+)_{\ast}(\alpha_1)-S_{\ast}(\alpha_2)$, where $\alpha_i$, $i=1,2$ are the curves defined in paragraph \ref{g}. Observe that the curve $\gamma-S_{\ast}(\gamma_2)$ bounds a disk in $\overline{\cal N}$ whose projection in the $z$-plane is the unit disk. So, we have
\begin{equation} \label{I} \int_{\gamma-S_{\ast}(\gamma_2)}\frac{dg}{g}= 2 \pi \ri \left({\rm Residue}\left(\frac{dg}{g},\ri a_+\right)+{\rm Residue}\left(\frac{dg}{g},-\ri b_+\right)\right)=2 \pi \ri \:(\beta+1) \, .
\end{equation}
Using again the notation of paragraph \ref{g} and taking into account (\ref{acciong}),  $\ri \frac{dg}{g}\left( \frac{d \delta_1}{dt} \right) \in \r$ and that $b$ and $a_3$ have been chosen to satisfy (\ref{eq3}) (or equivalently equations (\ref{eq1}) and (\ref{eq2})), we also obtain
\begin{equation} \label{II}\int_{\gamma-S_{\ast}(\gamma_2)}\frac{dg}{g}=\int_{\gamma+\gamma_2}\frac{dg}{g}=\int_{\gamma}\frac{dg}{g}
= 2 \int_{\delta-\delta_1-(S_2^+)_{\ast}(\delta-\delta_1)}\frac{dg}{g}=2 \int_{\delta-\delta_1} \frac{dg}{g}-\overline{\frac{dg}{g}}= 4  \int_{\delta-\delta_1} \frac{dg}{g}\; .
\end{equation}
Then, from (\ref{I}) and (\ref{II}) we get
$$ g(\ri_+)=\exp \left( \int_{\delta-\delta_1} \frac{dg}{g} \right)= \re^{ \frac{\ri (\beta+1)}{2}} \; .$$
Analogously, we can prove that $ g(-\ri_+)= \re^{ -\frac{\ri (\beta+1)}{2}}$. By using these computations and (\ref{Phi3}), it is not difficult to see that on the curves $\gamma_1^{\pm}$ we have
$$g(\gamma_1^{\pm}(t))=\re^{\pm \frac{\ri (\beta+1)}{2}} \psi(t) \; ,  \quad \Phi_3 \left( \frac{d\gamma_1^{\pm}}{dt}\right)=\ri \: \lambda \: \varphi(t) dt\; ,$$
where $\psi:[-1,-b[\: \cup \:]-b,a[ \: \cup \:]a,1] \longrightarrow \r$ and $\varphi:[-1,a[\: \cup \:]a,1] \longrightarrow \r$ have the following properties
\begin{itemize}
\item $\psi>0$ and $\varphi<0$ on $]a,1]$,
\item $\psi>0$ on $[-1,-b[$, $\psi<0$ on $]-b,a[$, $\varphi<0$ on $[-1,-b[$, $\varphi>0$ on $]-b,-a[$, $\varphi(b)=0$,
$\lim_{t \to -b^-} \psi=+ \infty$, $\lim_{t \to -b^+} \psi=- \infty$ and $\psi \varphi$ is well defined on $[-1,a[$ and  in this interval $\psi \varphi<0$.
\item $\lim_{t \to a^+} \varphi=-\infty$, $\lim_{t \to a^-} \varphi=\infty$, $\lim_{t \to a^+} \psi=\infty$ and $\lim_{t \to a^-} \psi=-\infty$.
\end{itemize}
Hence we have
$$\Phi_1\left( \frac{d\gamma_1^{\pm}}{dt}\right)=\frac{\ri \lambda}{2} \left( \re^{\mp \frac{\ri (\beta+1)}{2}} \frac{1}{\psi(t)}-\re^{\pm \frac{\ri (\beta+1)}{2}} \psi(t)\right) \varphi(t)\; , $$
$$\Phi_2\left( \frac{d\gamma_1^{\pm}}{dt}\right)=-\frac{\lambda}{2} \left( \re^{\mp \frac{\ri (\beta+1)}{2}} \frac{1}{\psi(t)}+\re^{\pm \frac{ \ri (\beta+1)}{2}} \psi(t)\right) \varphi(t)\; .$$
The expressions for $\Phi_i$, $i=1,2,3$, imply that $\ell_1^{\pm}$ is a half-line contained in a straight line  $\{ x_3=k^{\pm}, \: \: \mp \sin\left(\frac{\beta \pi}{2}\right)x_1+\cos\left(\frac{\beta \pi}{2}\right)x_2=K^{\pm} \}$, for suitable $k^{\pm}, K^{\pm} \in \r$. Observe that these straight lines meet the straight line $\{ x_3=k^{\pm}, \: \: x_1=0 \}$  at an angle $\frac{\beta \pi}{2}$. Note also that 
$${\rm Re}\left(\Phi_1\left( \frac{d\gamma_1^{\pm}}{dt}\right)\right)=\pm \frac{\lambda}{2}\sin\left(\frac{(\beta+1) \pi}{2}\right) \left( \frac{1}{\psi(t)}+\psi(t)\right) \varphi(t) \;. $$
$${\rm Re}\left(\Phi_2\left( \frac{d\gamma_1^{\pm}}{dt}\right)\right)=-\frac{\lambda}{2}\cos\left(\frac{(\beta+1) \pi}{2}\right) \left( \frac{1}{\psi(t)}+\psi(t)\right) \varphi(t) \;. $$ 
Taking into account the properties of the functions $\psi$ and $\varphi$ and that $\beta \in \:]0,1]$, we deduce that ${\rm Re}\left(\Phi_2\left( \frac{d\gamma_1^{\pm}}{dt}\right)\right)<0$ and so $X_{|\gamma_1^{\pm}}$ is injective. 

Moreover, it is clear that ${x_2}_{|\ell_1^{\pm}}$ diverges to $\pm \infty$. If $\beta \in \: ]0,1[$, then ${x_1}_{|\ell_1^{\pm}}$ diverges to $+\infty$, while $\beta=1$ implies that ${x_1}_{|\ell_1^{\pm}}$ is constant.

On the other hand, it is straightforward to check that $\Phi_3\left( \frac{d \delta_i}{dt} \right) \in \r$ and $\ri \frac{dg}{g}\left( \frac{d \delta_i}{dt} \right) \in \r$, for $i=1,2$. Moreover, since $b$ and $a_3$ had been selected to satisfy the system (\ref{eq3}) we have that the Weierstrass data verify equation (\ref{eq1}). All these facts imply that
$$\mbox{Re} \left( \Phi_j \left( \frac{d \delta _i}{dt} \right) \right) =0, \; i=1,2, \; j=1,2,$$
and so $\ell_0^+$ is a vertical segment.

Furthermore, note that
$$ \Phi_3 \left(\frac{d \delta _1}{dt} \right) \in \r^+, \; 
\Phi_3 \left( \frac{d \delta _2}{dt} \right) \in \r^-, $$
and this implies that $X_3(\ri_+)-X_3(\ri _-)>0$ and $X_{|\gamma_0^+}$ is injective. On the other hand, since $\ri_+ \in \gamma_1^+$ and $0 \in \gamma_1^-$ we have
$$ k^+-k^-=X_3(\ri_+)-X_3(0_+)= {\rm Re} \left( \int_{0_+}^{\ri_+} \Phi_3 \right)=\ri \: \pi \:{\rm Residue}(\Phi_3,\ri a_+)=\pi\: \lambda \: R \; ,$$
where $R=\frac{a^2}{1-a^4}\left( \sqrt{a^2+\frac{1}{a^2}+2 \cos(\rho)}+\sqrt{b^2+\frac{1}{b^2}+2 \cos(\rho)} \right) $. Thus,  $\ts =k^+-k^- \geq 0$. Taking into account that symmetry $S$ given in (\ref{newsim}) verifies $S(\gamma_1^+)=\gamma_2^-$, $S(\gamma_1^-)=\gamma_2^+$ and $S(\gamma_0^+)=\gamma_0^-$, it is not hard to conclude {\sc A.2} and that $X_{|\gamma^+}$ and $X_{|\gamma^-}$ are injective.

Next, we prove that $X$ is a proper immersion that satisfies {\sc A.4}. In order to do this we recall that $\frac{dg}{g}$ has a simple pole at $\ri a_+$ with residue $\beta$. Therefore, one has that the behavior of $g$ in a neighborhood of $\ri a_+$ in $M$ is given by
$$g=B_0(z-\ri a)^{\beta}+H_0(z) \, ,$$
where $B_0\neq 0$, $H_0$ is a holomorphic function in that neighborhood and the branch of $(z-\ri a)^{\beta}$ satisfies $1^{\beta}=1$. On the other hand, since $\Phi_3$ has a simple pole at $\ri a_+$ one has that in a neighborhood of $\ri a_+$ in $M$ 
$$\Phi_3(z)=\frac{-\ri \: \lambda \: R}{z-\ri a}+H_1(z) \; ,$$
where $H_1$ is a holomorphic function in that neighborhood. Hence we deduce that
$$\Phi_1(z)= \frac{B_1}{(z-\ri a)^{\beta+1}}+H_2(z) \; ,\quad \Phi_2(z)=  \frac{\ri B_1}{(z-\ri a)^{\beta+1}}+H_3(z) \; ,$$
where as before $H_2$ and $H_3$ are holomorphic functions at $\ri a_+$. From expressions of $\Phi_i$, for $i=1,2,3$ in a neighborhood of $\ri a_+$ we have
\begin{eqnarray*}X(z) & = & (X_1(z)+\ri X_2(z),X_3(z))\\
&= &\left({\rm Re}\left(\frac{B_1}{\beta(z-\ri a)^{\beta}}+H_4(z)\right)-
\ri \:{\rm Im}\left(\frac{B_1}{\beta(z-\ri a)^{\beta}}+H_5(z)\right),{\rm Re}\left( -\ri  \lambda  R \log(z- \ri a) +H_6(z) \right) \right)\\
&= &\left(\overline{\frac{B_1}{\beta(z-\ri a)^{\beta}}},0\right)+(O_1(1)+\ri O_2(1),O_3(1)) \; ,
\end{eqnarray*}
where $O_i(1)$ is a real function bounded in a neighborhood of $\ri a_+$ and $H_i$, for $i=4,5,6$ are holomorphic functions at $\ri a_+$. Firstly, from the above expression it follows that $X$ is proper. We also have that $\arg\left(\overline{\frac{B_1}{\beta(z-\ri a)^{\beta}}}\right) \in [\arg(\overline{B_1})-\frac{\pi}{2},\arg(\overline{B_1})+\frac{\pi}{2}]$ and then $X(M) \subset \Sigma^+$, where $\Sigma^+$ is a half space of $\r^3$ determinated by a plane $\Sigma$ orthogonal to $\{x_3=0\}$ that verifies $X(\partial(M)) \subset \Sigma^+$. Furthermore, we can infer from the preceding expression that $X(M)\subset \{(x_1,x_2,x_3)\: \mid \:k \leq x_3 \leq k'\}$, where $k< k'$. In the case $\beta=1$, we can use a result by Meeks and Rosenberg (see Lemma 2.1 in \cite{m-r}) to obtain that $X(M)$ lies in the half slab determinated by to horizontal planes containing $\ell_1^+$ and $\ell_2^-$, respectively, and $\Sigma_0^+$, where $\Sigma_0^+$ is the half space determinated by the plane parallel to $\Sigma$ that contains $\partial(M)$ and $\Sigma_0^+ \subset \Sigma^+$. This concludes the proof of {\sc A.4} in the case $\beta=1$. If $\beta \in ]0,1[$ we have that $X(M)$ is in a half slab and $X(\partial(M))$ is contained in a wedge of angle $\pi \beta \in \:]0,\pi[$. Then a result by López and Martín (see Corollary 2 in \cite{l-m1}) asserts that with this conditions $X(M)$ lies in the convex hull of its boundary. 
\end{proof}
\begin{remark} \label{paco} Observe that the case $a=0$ in our family corresponds with the López-Martín examples (see \cite{l-m1}). 
\end{remark} 

The second objective of the present section is to prove the following result:
\begin{theorem} \label{gordo}
For each $\beta \in ]0,1]$ there exists $(a(\beta), \rho(\beta))$, such that 
\begin{equation} \sh (a(\beta),\rho(\beta), \beta)=0 \; , \quad \sd(a(\beta), \rho(\beta), \beta)=0  \; .\label{eqdh}
\end{equation}
\end{theorem}

Assume $\beta \in ]0,1]$ fixed. First we try to write  the system (\ref{eqdh}) in terms of integrals of the Weierstrass data. Observe that $0_+ \in \gamma_1^-$ and $\infty_- \in \gamma_2^+$. On the other hand we have the symmetry $S$ that satisfies $S(0_+)=S(\infty_-)$ and $S(1_+)= 1_+$. Thus, taking into account (\ref{simphi3}) we get the following expression for the function $\sh$
$$ \sh(a,\rho, \beta)=x_3(q_2^+)-x_3(q_1^-)=X_3(\infty_-)-X_3(0_+)={\rm Re}\left(\int _{0_+}^{\infty_-} \Phi_3 \right)=2 \:{\rm Re}\left( \int _{0_+}^{1_+} \Phi_3 \right)\,.  $$
So the first equation in (\ref{eqdh}) is equivalent to 
\begin{equation} \label{eqh}
h(a,\rho, \beta)={\rm Re}\left( \int _{0_+}^{1_+} \frac{1}{\lambda} \Phi_3 \right)=0\,.  
\end{equation}
Regarding the function $\sd$, it is clear that $\re^{\ri \frac{\rho}{2}} \in \gamma_0^+$ and $\re^{-\ri \frac{\rho}{2}} \in \gamma_0^-$, thus we have
\begin{equation}\label{eqsd}
 \sd(a,\rho, \beta)=X_2(\re^{\ri \frac{\rho}{2}})-X_2(\re^{-\ri \frac{\rho}{2}})=-\frac{1}{2}\: {\rm Im} \left( \int _{\re^{-\ri \frac{\rho}{2}}}^{\re^{\ri \frac{\rho}{2}}} \left( g+\frac{1}{g} \right) \Phi_3 \right)= -\frac{1}{2}\: {\rm Im} \left( \int _{\alpha_1} \left( g+\frac{1}{g} \right) \Phi_3 \right)\; , 
\end{equation}
where $\alpha_1$ is the curve defined in paragraph \ref{g}. We recall that $\alpha_1=-S_{\ast}(\delta)+\delta$ and thus $S_{\ast}(\alpha_1)=-\alpha_1$. Then, taking into account (\ref{simphi3}) and (\ref{acciong}) we deduce that \begin{equation}\label{sinponer}
 \int_{\alpha_1} \frac{1}{g} \Phi_3=\int_{\alpha_1} g \Phi_3 \;.
\end{equation}
From here and (\ref{eqsd}) 
we obtain that the second equation in (\ref{eqdh}) is equivalent to 
\begin{equation} \label{eqd}
d(a,\rho, \beta)=- \frac{1}{\lambda}{\rm Im} \left( \int _{\alpha_1} g \Phi_3 \right)=- \frac{1}{\lambda}{\rm Im} \left( \int _{\re^{-\ri \frac{\rho}{2}}}^{\re^{\ri \frac{\rho}{2}}} g \Phi_3 \right)=0 \; .
\end{equation}
Thus, to prove Theorem \ref{gordo} it suffices to prove that there exists $(a(\beta), \rho(\beta))$ such that
$$ h(a(\beta),\rho(\beta), \beta)=0 \; , \quad d(a(\beta), \rho(\beta), \beta)=0  \, .$$
Unfortunately, the proof of the above assertion is quite long and technical. For the sake of clarity, we develop here a sketch of the proof and we present the complete details in Sect. \ref{detalles}.

We study first the function $h$. From (\ref{eqh}) we have that
\begin{equation}\label{hdef} h(a,\rho,\beta)= \int_0^1 \frac{t^4 +1+(c_1 c_2-2 \cos \rho) t^2}{(t^2+a^2) (t^2+\frac{1}{a^2})\sqrt{t^4+1-2 t^2 \cos \rho}} dt \; . 
\end{equation}
Observe that the function $h$ can be extended to $[0,1] \times ]0,\pi] \times ]0,1]$. We obtain the following properties for $h$:
\begin{enumerate}[{\sc {Claim} H.1}] 
\item $h(0,\rho,\beta) <0$ for $\rho \in \: ]0,\pi]\; ,$ 
\item $\lim_{\rho \to 0^+}h(a,\rho,\beta)= -\infty$ for $a \in [0,1] \; ,$ 
\item $h(1,\pi,\beta)=\frac{\pi}{4} \; ,$
\item $\frac{\partial h}{\partial \rho}(a,\rho,\beta) > 0 $  for $(a,\rho,\beta) \in \:]0,1] \times ]0,\pi[ \times ]0,1]  \; ,$
\item $\frac{\partial h}{\partial a}(a,\rho,\beta) > 0$ for $(a,\rho,\beta) \in \:]0,1[ \times ]0,\pi] \times ]0,1]  \; .$
\end{enumerate}
From the preceding assertions we have that there exist $\rho_0(\beta)\in \:]0,\pi[$ and $a_0(\beta)\in \:  ]0,1[$  such that $h(1,\rho_0(\beta),\beta)=h(a_0(\beta),\pi,\beta)=0$. Furthermore, the set $ C_1=\{ (a,\rho) \in [0,1] \times [0,\pi] \: \mid \: h(a ,\rho,\beta)=0 \}$ is a differentiable embedded curve in $]0,1] \times ]0,\pi]$ from $(a_0(\beta), \pi)$ to  
$(1,\rho_0(\beta))$. 

However, the study of the function $d$ is much more complicated due to the expression of the Gauss map $g$. For this function the following facts can be proved:

\begin{enumerate}[{\sc {Claim} D.1}] 
\item $d(a,0,\beta) > 0 \; ,$ 
\item $\lim_{\rho \to \pi}d(a,\rho,\beta) = -\infty \; ,$
\item There exists a unique $\rho_1(\beta) \in ]0, \pi[$ such that $d(0,\rho_1(\beta),\beta)=0$. Furthermore, $d$ is positive for $\rho \in [0,\rho_1(\beta)[$ and negative for $\rho \in \: ]\rho_1(\beta),1]\; ,$
\item $d(1,\rho,\beta)>0$ for $\rho \in [0,\rho_0(\beta)]\; .$
\end{enumerate}
All the preceding claims allow us to assert that there exists a connected component $C$ of the set $ \{ (a,\rho) \in [0,1] \times [0,\pi] \: \mid \: d(a ,\rho,\beta)=0 \}$ that contains the point $(0,\rho_1(\beta))$ and a point in the segment $\{(1,\rho) \: \mid \: \rho \in [\rho_0(\beta),\pi] \}$.  Since $d$ is a real analytic function  in the interior of $[0,1] \times [0,\pi]$, we have that $d^{-1}(0)$ is locally path-connected and thus $C$ is also a path component. We denote by $C_2$ a path in $C$ starting at $(0,\rho_1(\beta))$ and finishing at a point in the segment $\{(1,\rho) \: \mid \: \rho \in [\rho_0(\beta),\pi] \}$ . Therefore, $C_1 \cap C_2 \neq \emptyset$ and this concludes the proof of Theorem \ref{gordo}.

\begin{figure}[htbp]
	\begin{center}
		\includegraphics[width=0.50\textwidth]{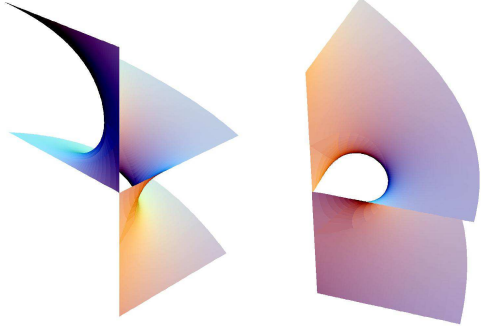}
	\end{center}
	\caption{A fundamental piece ($\sd=\sh=0$) of angle $\pi/2$.}
	\label{fig:trozo}
\end{figure}

 \begin{theorem} \label{embebi}
If $(a, \rho) \in C_2$, then  the surface given by the Weierstrass data (\ref{Phi3}) also fulfills {\sc A.5}.
\end{theorem}
\begin{proof}
Let $C_2(s)=(a(s),\rho(s))$, $s \in [0,1]$ a parametrization of the curve $C_2$. We label $M_s=M_{a(s),\rho(s)}$ and $X_s=X_{a(s),\rho(s),\beta}$. We also define the set
${\cal I} =\{s \in [0,1] \: \mid   \: X_s \quad \mbox{satisfies {\sc A.5}} \} \; .$ To conclude it is sufficient to see that ${\cal I}=[0,1]$.

Recall that $a=0$ corresponds to the López-Martín example of angle $\pi \beta$. This surface verifies {\sc A.5} when $d=0$ (see \cite{l-m1}). So, $0 \in {\cal I} \neq \emptyset$.
Consider $s_0 \in {\cal I}$. Taking into account the Weierstrass representation,  one has that the ends of $X_{s}$ are asymptotic to two pieces of helicoid and the distance between both ends is positive, $\forall s$.  Then, there exists $\epsilon,r>0$ such that $X_s (M) \cap (\r^3-B(0,r))$ is  embedded, $\forall s \in ]s_0-\epsilon, s_0+\epsilon[$. If $X_s$ were not injective for some $s \in  ]s_0-\epsilon, s_0+\epsilon[$, then self-intersections of $X_s(M_s)$ would be in $ \overline{B}(0,r)$. As  $X_s(M_s)$ is contained in the convex hull of its boundary, there are no contacts of interior points with points at the boundary. So we would arrive to a contradiction, by using either the classical maximum principle or the maximum principle at the boundary. Hence, $]s_0-\epsilon, s_0+\epsilon[ \subset {\cal I}$ which implies that $\cal I$ is open.

Now, take $\{s_n\}_{n \in \n}$ a sequence in $\cal I$ converging to $s_0>0$. Assume that $X_{s_0}$ is not injective. Then, there are two points $x,y \in M_{s_0}$ satisfying $X_{s_0}(x)=X_{s_0}(y)$. The convergence of $\{X_{s_n}\}_{n \in \n}$ to $X_{s_0}$ uniformly over compact subsets of $\r^3$ and the interior maximum principle assure that there exist neighborhoods
$N(x)$, $N(y)$, of $x$ and $y$, respectively, such that $X_{s_0}(N(x))=X_{s_0}(N(y))$. So, the image set $X_{s_0}(M)$ is an embedded minimal surface with finite total curvature and 
$X_{s_0}:M_{s_0} \longrightarrow X_{s_0}(M_{s_0}) $
is a finitely sheeted covering map. As $X_{s_0}$ is one-to-one in a neighborhood of the end, then we deduce that $X_{s_0}$ is injective, which is contrary to our assumption. This contradiction proves that $\cal I$ is closed.

Thus, an elementary connectedness argument gives that ${\cal I}=[0 , 1]$, which concludes the proof.
\end{proof}
We can describe our family of surfaces $\cal M$ as follows:
$$ {\cal M}=\{X_{a,\rho, \beta}: M_{a,\rho} \longrightarrow \r^3 \: \mid \: a \in [0,1[, \rho \in \: ]0,\pi[, \beta \in \:]0,1] \} \; .$$

Finally, we are interested in the complete surfaces obtained from the minimal immersion $X_{a,\rho, \beta}: M_{a,\rho} \longrightarrow \r^3$ when the parameters $(a,\rho,\beta)$ satisfy (\ref{eqdh}). We summarize the properties of these complete surfaces in the following remark.
\begin{remark} \label{re:singly}
As we mentioned in the paragraph \ref{M}, the complete orientable  minimal surface without boundary 
$$\widetilde{X}_{a,\rho,\beta}:\widetilde{M}_{a,\rho,\beta} \longrightarrow \r^3$$ 
obtained  from $X_{a,\rho,\beta} (M_{a,\rho})$
 by successive Schwarz reflections about straight lines is invariant under the vertical translation $\tau_1$. 

The case $\beta \in \q$ is specially interesting. 
The immersion $\widetilde{X}_{a,\rho,\beta}$ is singly periodic and the induced
immersion 
$$Y_{a,\rho,\beta}:(\widetilde{M}_{a,\rho,\beta})/\langle \tau_1 \rangle \longrightarrow \r^3/\langle \tau_1 \rangle$$
has four ends and finite total curvature.  If we write $\beta=q/p$, $p, \;q \in \n$, $gcd(p,q)=1$, then
it is not hard to check that:  
\begin{itemize}
\item If $p$ is even the surface $\widetilde{M}_{a,\rho,\beta}$ is the two sheeted orientable covering 
of a nonorientable minimal surface properly immersed in $\r^3$.  Moreover, $Y_{a,\rho,\beta}$ has four ends, 
its total curvature is $-8 \pi (p+q)$ and $\widetilde{M}_{a,\rho,\beta}/\langle \tau_1 \rangle$ has genus $2p-1$.
A fundamental piece bounded by straight lines of a surface $\widetilde{X}_{a,\rho,\beta}$,
$\beta=1/2$, is illustrated in Figure \ref{fig:heli-2}.

\item If $p$ is odd  the induced immersion 
$Y_{a,\rho,\beta}$ has two ends. Moreover,  if $q$ is even (resp. $q$ is odd),  $Y_{a,\rho,\beta}$ has 
total curvature  $-8 \pi (p+q)$  (resp. $-4 \pi (p+q)$) and $\widetilde{M}_{a,\rho,\beta}/\langle \tau_1 \rangle$ 
has genus $2p$ (resp. $p$). Figures \ref{fig:heli-3} and \ref{fig:heli-23} shows a fundamental piece of a surface
$\widetilde{X}_{a,\rho,\beta}(\widetilde{M}_{a,\rho,\beta})$, $\beta=1/3$ and $\beta=2/3$, respectively.
\begin{figure}[h]
	\begin{center}
		\includegraphics[width=0.50\textwidth]{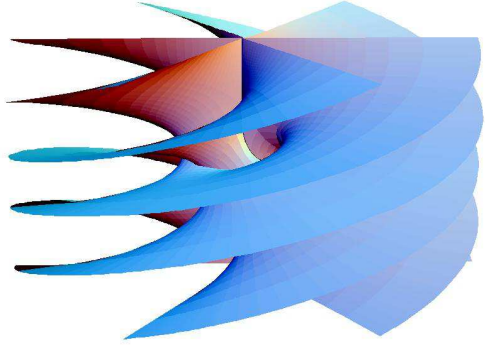}
	\end{center}
	\caption{A complete surface constructed with a fundamental piece of angle $2 \pi/3$.}
	\label{fig:heli-23}
\end{figure}
\item The ends of $Y_{a,\rho,\beta}$ are embedded if, and only if, $\beta=\frac{1}{n}$, $n \in \n$, $n>1$. 
In this case the only self-intersection of $\widetilde{X}_{a,\rho,\beta}(\widetilde{M}_{a,\rho,\beta})$ occurs along the $x_3-$axis. 
\end{itemize}
\end{remark}

\begin{remark}[Uniqueness of ${\cal H}_1$] \label{hkw} Observe that the case $\beta=1$ (i.e., angle $\pi$) corresponds to the simply periodic genus-one helicoid whose existence and embeddedness were proved by Hoffman-Karcher-Wei in \cite{hkw}.
 
From (\ref{F}), it is clear that if $\beta=1$  we have $b(a,\rho,1)=a$. Hence $a_2=-a_1$, $c_2=-c_1$ and $a_3=0$. Thus we obtain the following expressions for the Weierstrass data in (\ref{Phi3})
$$ g=\frac{z^2+a^2}{a^2 z^2+1} \; , \quad \Phi_3= \lambda \frac{w-c_1 \ri z}{w+c_1 \ri z}\frac{dz}{w} \; .$$
It is easy to check that 
$$ \frac{1}{\lambda} g \Phi_3= d\left( \frac{a^2 z w+\ri c_1}{a^2 z^2+1}\right)-\frac{z^2}{w} dz \; .$$
Thus, from the expression (\ref{dgamma1}), that we will prove in the following section, we have
$$d(a,\rho,1)=-\frac{\ri}{2} \int_{\gamma_1}\frac{z^2}{w} dz \; .$$
Therefore, the function $d$ in this case does not depend on $a$. Then, taking into account {\sc Claim H.3}, the set $ \{ (a,\rho) \in [0,1] \times [0,\pi] \: \mid \: d(a ,\rho,1)=0 \}=\{ (a,\rho_1(1)) \: \mid \: a \in [0,1] \}=C_2$. Finally, using {\sc Claim H.5} and that $C_1$ is a graph in $\rho$ we obtain that $C_1 \cap C_2$ is a unique point and so we have the uniqueness of the example when the angle is $\pi$.
\end{remark} 
\begin{remark}[Limit case $\beta=0$] \label{js} Regarding to the case $\beta=0$, we observe that $\lim_{\beta \to 0} b(a,\rho,\beta)=0$ and therefore $\lim_{\beta \to 0} a_3(a,\rho,\beta)=0$ (see Lemma \ref{lemab} and (\ref{a31})). Furthermore, if we impose, in order to take limits, that  the length of the vertical segments were $1$ we obtain that 
$$\lambda= {\rm Re}\left(\int_{\infty_-}^{\ri_+}  \frac{w+c_2 \ri z}{w+c_1 \ri z} \frac{dz}{w}\right)
=-2h(a,\rho,\beta)+{\rm Re}\left(\int_{0}^{\ri_+}  \frac{w+c_2 \ri z}{w+c_1 \ri z} \frac{dz}{w}\right)=-2h(a,\rho,\beta)+ \pi R \; ,$$
where $R$ was defined in the proof of Theorem \ref{gordo}. Hence, if we study the limits of the Weierstrass data given in (\ref{Phi3}) as $\beta \to 0$ we obtain
$$g=z \; ,\;
\Phi_3= \mu \frac{\ri z}{w+c_1 \ri z} \frac{dz}{w} \; , $$
where $\mu \in \r$. One can see, using similar arguments as in Theorem \ref{teo:grande}, that the minimal surface with  this Weierstrass data satisfies assumptions {\sc A.1} to {\sc A.4}. Taking into account the expression of the Gauss map, it is easy to prove that this surface is a graph over the plane $\{x_1=0\}$. Therefore, this Weierstrass data corresponds to a Jenkins-Serrin graph defined on a rhomboid (see Fig. \ref{fig:romboide}).
\begin{figure}
	\begin{center}
		\includegraphics[width=0.15\textwidth]{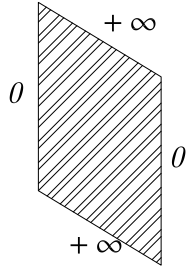}
	\end{center}
	\caption{The rhomboid with the boundary values for the Jenkins-Serrin graph.}
	\label{fig:romboide}
\end{figure}

\end{remark}
\section{Appendix: Technical computations}\label{detalles}
The aim of this section is to present in detail the proofs of the claims in Sect. \ref{completas}. Firstly, we have to study in deep the function $b(a,\rho,\beta)$ defined in the paragraph \ref{g}. More precisely, we prove:
\begin{lemma} \label{lemab} The function $b:[0,1] \times \: ]0,\pi] \times [0,1] \longrightarrow [0,a]$ fulfills the following properties:
\begin{enumerate}[{a)}]
\item $b(a,\rho,\beta)=0$ if and only if $a=0$ or $\beta=0 \, ,$
\item $b(a,0,\beta)=\tan(\beta {\rm arctan}(a))$ and $b(a,\pi,\beta)=\tanh(\beta {\rm arctanh}(a)) \, ,$
\item $\lim _{a \to 0}\frac{b(a,\rho,\beta)}{a}=\beta \, ,$
\item $\frac{\partial b}{\partial \rho}(a,\rho,\beta) \geq 0$ for $(a,\rho,\beta) \in \:]0,1] \times ]0,\pi[ \times [0,1] \; , $
\end{enumerate}
\end{lemma}
\begin{proof} The assertions \emph{a)} and \emph{b)} are obtained straightforward from (\ref{F}). Moreover, taking into account that expression (\ref{F1}) vanishes for $b=b(a,\rho,\beta)$, it is not difficult to see \emph{c)}. Now we prove assertion \emph{d)}.

For the sake of brevity, we will denote $b_{\rho}=\frac{\partial b}{\partial \rho}$, $b_{\beta}=\frac{\partial b}{\partial \beta}$, $b_{\rho \beta}=\frac{\partial^2 b}{\partial \rho \partial \beta}$ and $b_{\rho \beta \beta}=\frac{\partial^3 b}{\partial \rho \partial \beta \partial \beta}$. By deriving in (\ref{F}) we obtain the following expressions
\begin{eqnarray}\label{br}
 b_\rho &=& \sin \rho \sqrt{b^4+1+2 b^2\cos \rho }\left(\beta \int_0^a  \frac{t^2 dt}{(t^4+1+2t^2 \cos \rho )^{\frac{3}{2}}} -\int_0^b  \frac{t^2 dt}{(t^4+1+2t^2 \cos \rho )^{\frac{3}{2}}}\right)\; ,\\
b_\beta &=& \sqrt{b^4+1+2 b^2\cos \rho } \int_0^a  \frac{dt}{t^4+1+2t^2 \cos \rho } \; . \label{bb}
\end{eqnarray}
Deriving again in the preceding expressions we get the following differential equations
\begin{eqnarray}
b_{\rho \beta}&=& A+ B \, b_{\rho} \, , \label{brb} \\
b_{\rho \beta \beta}&=& C \,b_{\beta}\, b_{\rho \beta}+D \,b_{\beta}^2 \, b_{\rho}+E \, b_{\beta}^2 \; , \label{brbb}
\end{eqnarray}
where
\begin{eqnarray*} A &=& \sin \rho \sqrt{b^4+1+2 b^2\cos \rho } \int _{0}^a \frac{t^2 dt}{\sqrt{t^4+1+2t^2 \cos \rho }}
- \frac{b^2 \sin \rho}{ \sqrt{b^4+1+2 b^2\cos \rho }} \int_0^a \frac{dt}{\sqrt{t^4+1+2 t^2\cos \rho }} \; ,\\
B &=& \frac{2\, b \,(b^2+\cos \rho)}{ \sqrt{b^4+1+2b^2 \cos \rho }} \int_0^a \frac{dt}{\sqrt{t^4+1+2t^2 \cos \rho }} \; ,
 \quad \quad C =  \frac{4 b (b^2+ \cos \rho)}{b^4+1+2 b^2\cos \rho } \; ,\\
D &=&  \frac{2\left( 1 - b^4 \right) \,(b^2+\cos (\rho ) )+4 \, b^2\,{\sin (\rho )}^2 }{(b^4+1+2 b^2\cos \rho)^2 }\; , \quad \quad E = -\frac{2 \, b\,{\sin (\rho )\,\left( 1 - b^4 \right)} }{(b^4+1+2 b^2\cos \rho)^2 }\; . 
\end{eqnarray*}
Let us prove that there does not exist any point $(a,\rho,\beta) \in \:]0,1] \times ]0,\pi[ \times [0,1]$
 such that 
\begin{equation} b_{\rho}(a,\rho, \beta)<0 \;, \quad b_{\rho \beta}(a,\rho, \beta)=0 \quad {\rm and} \quad b_{\rho \beta \beta}(a,\rho, \beta)\geq 0 \; .\label{tres}
\end{equation}
Assume we have a point $(a,\rho, \beta)$ with $b_{\rho}(a,\rho, \beta)<0$ and $b_{\rho \beta}(a,\rho, \beta)=0$. If $ \rho \in ]0, \frac{\pi}{2}]$ we obtain from (\ref{bb}), (\ref{brbb}) and the definitions of $D$ and $E$ that $b_{\rho \beta \beta}(a,\rho, \beta) < 0$. 

Suppose now that $\rho \in [\frac{\pi}{2},\pi]$. From our assumptions and (\ref{brb}) we have  
\begin{equation}b(a,\rho, \beta) \sin \rho \geq 2 \:b_\rho(a,\rho, \beta) (b(a,\rho, \beta)^2+ \cos \rho)\; .\label{ine}
\end{equation}
Then, taking into account (\ref{brbb}), (\ref{ine}) and that $b \in [0,1]$  it is clear that 
$$b_{\rho \beta \beta}(a,\rho, \beta)  \leq  \left(\frac{2\:(b^2-\cos \rho)}{b^4+1+2 b^2\cos \rho }b_\beta^2 \: b_\rho \right)(a,\rho, \beta)<0 \; .$$
On the other hand, from (\ref{br}) we have $b_\rho(a,\rho,0)=b_\rho(a,\rho,1)=0$. Since there are no points satisfying (\ref{tres}) we deduce that $b_\rho$ as a function of $\beta$ have no minimums corresponding to negative values of $b_\rho$. Hence $b\rho \geq 0$ and this concludes the proof.
\end{proof}
{\sc Proof of claim H.1 .} This assertion can be inferred straightforward from (\ref{ces}), (\ref{hdef}) and parts \emph{a)} and \emph{c)} of Lemma \ref{lemab}.\\
\noindent{\sc Proof of claim H.2 .} Note that if $\rho=0$, the denominator of (\ref{hdef}) has a zero of order one at $t=1$ while the numerator is strictly negative at this point. This gives the claim.

\noindent{\sc Proof of claim H.3 .} From part \emph{a)} in Lemma \ref{lemab} and (\ref{ces}) we have that $c_1 \cdot c_2(1,\pi,\beta)=0$. Therefore the integral (\ref{hdef}) becomes
$$h(1,\pi,\beta)=\int_0^1 \frac{dt}{t^2+1}=\frac{\pi}{4} \; .$$
\noindent{\sc Proof of claim H.4 .} To prove this claim we will consider the function $\tilde{h}(a,b,\rho)= h(a,\rho,\beta)$, it is to say we consider $\tilde{h}$ as a function of the independent variables $(a,b,\rho)$. Then
$$\frac{\partial h}{\partial \rho}=\frac{\partial \tilde{h}}{\partial \rho}+\frac{\partial \tilde{h}}{\partial b} b_\rho \, .$$
By deriving in (\ref{hdef}) we get
$$ \frac{\partial \tilde{h}}{\partial \rho}= \frac{-\sin \rho }{c_1 c_2}\int_0^1 \frac{t^2 \,\left( (c_1^2+c_2^2-c_1 c_2)(t^4+1-2 t^2 \cos \rho)+c_1^2 c_2^2 t^2 \right)}{(t^2+a^2)(t^2+\frac{1}{a^2}) (t^4+1-2 t^2 \cos \rho)^{\frac{3}{2}}}dt \; ,$$
$$ \frac{\partial \tilde{h}}{\partial b}= \frac{- c_1 (1-b^4)}{ c_2\: b^3} \int_0^1 \frac{t^2}{(t^2+a^2)(t^2+\frac{1}{a^2}) \sqrt{t^4+1-2 t^2 \cos \rho}}dt \; .$$
From here and part \emph{a)} in Lemma \ref{lemab} we have $\frac{\partial \tilde{h}}{\partial \rho}>0$ and 
$\frac{\partial \tilde{h}}{\partial b}\geq 0$ for $(a,\rho,\beta) \in \: ]0,1] \times ]0,\pi[ \times ]0,1]$ .  Therefore, taking into account part \emph{d)} in Lemma \ref{lemab} we obtain {\sc Claim H.4}.

\noindent{\sc Proof of claim H.5 .} By deriving in (\ref{hdef}) and simplifying we obtain
$$\frac{\partial h}{\partial a}=  \left(\frac{1}{a^3}-a \right) \int_0^1 \frac{t^2 \left( 2-\frac{c_2}{c_1} t^2 \right)\sqrt{t^4+1-2 t^2 \cos \rho}}{(t^2+a^2)^2(t^2+\frac{1}{a^2})^2 }dt +\int_0^1 \frac{t^2 P(t,a,\rho,\beta)}{(t^2+a^2)^2(t^2+\frac{1}{a^2})^2 \sqrt{t^4+1-2 t^2 \cos \rho}}dt \; , $$
where $P(t,a,\rho,\beta)=A t^4 + B t^2+ A$ and $A= \frac{\beta b}{a} \left(\frac{1}{b^3}-b \right)$, $B=c_1 c_2 \left( \frac{1}{a^3}-a \right)+A \left(\frac{1}{a^2}+a^2 \right)$. Evidently, the first summand in the above expression is positive for $ a \in \: ]0,1[$. Then it suffices to see $P(t,a,\rho,\beta) \geq 0$. Since $A \geq 0$, if $B \geq 0$ this is obvious. Therefore, assume $B < 0$. In order to see $P(t,a,\rho,\beta) \geq 0$ we observe that the discriminant of $P$ as a polynomial in $t$ is given by $\Delta =B^2-4 A^2= (B-2 A)(B+2 A)$. From our assumption we know that the first factor of $\Delta$ is non positive. Let us analyze the second one. We have
$$B+2 A=c_1 c_2 \left( \frac{1}{a^3}-a \right)+\frac{\beta b}{a} \left(\frac{1}{b^3}-b \right)\left(\frac{1}{a}+a \right)^2 \; .$$
Taking into account that $x^4+1-2 \: x^2 \cos \rho \leq (x^2+1)^2$ we obtain
$$B+2 A \geq \frac{(1+a^2)^2 (1+b^2)}{a^3 b}\left(\beta \: \frac{1-b^2}{b} -\frac{1-a^2}{a} \right) \, .$$
Let us consider the function $f_1(a,\rho,\beta)= \beta \: \frac{1-b^2}{b} -\frac{1-a^2}{a}$. From parts \emph{b)} and \emph{d)} in Lemma \ref{lemab} we have that $b(a,\rho,\beta) \leq b(a,\pi,\beta)= \tanh(\beta \: {\rm arctanh}(a))$. Since the function $\frac{1-b^2}{b}$ is decreasing in $b$ we get
$$ f_1(a,\rho,\beta) \geq \beta \frac{2}{\sinh (2 \:\beta \: {\rm arctanh}(a))}-\frac{1-a^2}{a} \; .$$
Next we compute
$$\frac{\partial f_1}{\partial \beta}= \frac{2 \cosh (2 \:\beta \: {\rm arctanh}(a))}{\sinh (2 \:\beta \: {\rm arctanh}(a))^2}(\tanh (2 \:\beta \: {\rm arctanh}(a))-2 \:\beta \: {\rm arctanh}(a)) \; . $$
Since $\tanh(x)-x \leq 0$ we deduce that $\frac{\partial f_1}{\partial \beta} \leq 0$ and so 
$$ f_1(a,\rho,\beta) \geq  f_1(a,\rho,1)=  \frac{2}{\sinh (2 \: {\rm arctanh}(a))}-\frac{1-a^2}{a}=0 \; .$$
Thus $B+2 A \geq 0$ and then $\Delta \leq 0$. As $P(0,a,\rho,\beta)=A \geq 0$, this implies that $P(t,a,\rho,\beta)\geq 0$ and this concludes the proof of the claim.

Next we prove assertions on function $d$. In order to do this we need to establish some previous results. The first one consists of finding an upper bound for the point $\rho_0(\beta)$ defined in Sect. \ref{completas}. To be more precise we can see
\begin{lemma}\label{lemaro} $\rho_0(\beta) \leq \frac{\pi}{\beta+1}$.
\end{lemma}
\begin{proof}Taking into account {\sc Claim H.5}, it suffices to prove that $H(\beta)=h(1,\frac{\pi}{\beta+1},\beta)\geq 0$. To prove this fact we will show that 
$$ -2\cos \left(  \frac{\pi}{2(\beta+1)}\right) \sqrt{b^2+ \frac{1}{b^2}+2 \cos \left( \frac{\pi}{\beta+1} \right)}-2\cos \left( \frac{\pi}{\beta+1} \right) \geq -2 \; ,$$ 
or equivalently
$$ 2 \sin \left(\frac{\pi}{2(\beta+1)}\right)^2 -\cos \left(\frac{\pi}{2(\beta+1)} \right)\sqrt{b^2+ \frac{1}{b^2}+2 \cos \left(\frac{\pi}{\beta+1} \right)} \geq 0 \; .$$ 
Thereby, consider the function
$$ H_1(\beta)=4 \sin\left(\frac{\pi}{2(\beta+1)}\right)^4 -\cos \left(\frac{\pi}{2(\beta+1)}\right)^2 \left(b^2+ \frac{1}{b^2}+2 \cos\left(\frac{\pi}{\beta+1}\right)\right) $$

Using parts \emph{b)} and \emph{d)} of Lemma \ref{lemab} we obtain that $\tan (\frac{\beta \pi}{4}) \leq b(1,\frac{\pi}{\beta+1},\beta) \leq 1$. So, as the function $x^2+\frac{1}{x^2}$ is decreasing in $]0,1]$ we get after some computations
$$ H_1(\beta) \geq H_2(\beta)= 4 \left( \tan \left(\frac{\pi}{2(\beta+1)} \right)^2 -\sin \left( \frac{\beta \pi}{2} \right)^{-2} \right) \; . $$
By deriving in the above expression we have
$$ \frac{\partial H_2}{\partial \beta}(\beta)= 4 \pi \left( \frac{ \cos \left(\frac{\beta\pi}{2} \right)}{  \sin \left(\frac{\beta\pi}{2} \right)^3} -\frac{\frac{1}{(\beta+1)^2} \sin \left(\frac{\pi}{2(\beta+1)} \right)}{  \cos \left(\frac{\pi}{2(\beta+1)} \right)^3}    \right)=\frac{4 \pi}{\beta^2} \left( \frac{\beta^2 \cos \left(\frac{\beta\pi}{2} \right)}{  \sin \left(\frac{\beta\pi}{2} \right)^3} -\frac{\frac{\beta^2}{(\beta+1)^2} \cos \left(\frac{\beta\pi}{2(\beta+1)} \right)}{  \sin \left(\frac{\beta \pi}{2(\beta+1)} \right)^3}    \right)\;. $$
Now we define the function
$H_3(x)=\frac{x^2 \cos \left(\frac{\pi}{2}x \right)}{  \sin \left(\frac{\pi}{2}x \right)^3}$ for $x \in [0,1]$. Let us see that $H_2$ is a decreasing function in $x$. Observe that 
$$\frac{\partial H_3}{\partial x}(x)=\frac{x}{2}  \csc \left(\frac{\pi}{2}x \right)^4 \left(2 \sin(\pi x)-\pi x (2+ \cos(\pi x))\right) \; .$$
It is easy to see that $H_4(x)=2 \sin(\pi x)-\pi x (2+ \cos(\pi x))$ is a concave function with $H_4(0)=0$ and $H_4(1)=-\pi$. Therefore, we deduce that $H_3$ is decreasing. Hence, since $\frac{\beta}{\beta+1} \leq \beta$ we obtain
$$ \frac{\partial H_2}{\partial \beta}(\beta)= \frac{4 \pi}{\beta^2} \left( H_3(\beta)-H_3 \left(\frac{\beta}{\beta+1} \right) \right) \leq 0 \; .$$
Finally, taking into account $H_2(1)=0$ we conclude that $H_1(\beta) \geq H_2(\beta) \geq 0$.
\end{proof}
Next, we will study in deep the function $a_3(a,\rho,\beta)$ presented in paragraph \ref{g}. 
\begin{lemma} \label{lemaa} The function $a_3:[0,1] \times ]0,\pi] \times [0,1] \longrightarrow \r$ fulfills the following properties:
\begin{enumerate}[{a)}]
\item $a_3(a,0,\beta)=0\, ,$
\item $a_3(a,\pi,\beta)=0\, ,$
\item $a_3(1,\rho,\beta)=-\frac{1}{a_2} A_0$, where $0 \leq A_0 \leq 1 \; .$ 
\end{enumerate}
\end{lemma}
\begin{proof} Recall that $a_3$ can be computed either by the expression (\ref{a31}) or expression (\ref{a32}). Observe that the denominator of (\ref{a31}) diverges to $\infty$ as $\rho \to \pi$ while the numerator goes to a constant. Therefore we obtain part \emph{b)} of the lemma. On the other hand, when $a \to 1$ the expression (\ref{a31}) gives us $a_3(1,\rho,\beta)=-\frac{1}{a_2} A_0$, where
$$A_0= \frac{\int_0^\rho \frac{\sqrt{2(\cos(t)-\cos(\rho))}}{b^2+\frac{1}{b^2}+2 \cos(t)}dt}{\int_0^{\rho} \frac{1}{\sqrt{2(\cos(t)-\cos(\rho))}}dt} $$
and $b=b(1,\rho,\beta)$. Clearly $A_0 \geq 0$. Furthermore, checking that the denominator in the above fraction is greater than the numerator we have that $A_0 \leq 1$ and so part \emph{c)} is proved. 

As before the denominator in the expression (\ref{a32}) diverges to $\infty$ as $\rho \to 0$ while the numerator goes to a constant. Thus we get statement \emph{a)} of the lemma.
\end{proof}
\noindent{\sc Proof of claim D.1 .} From (\ref{eqd}) we have
$$ d(a,\rho,\beta)= - \frac{1}{\lambda}{\rm Im} \left( \int _{\alpha_1} g \Phi_3 \right) \; ,$$
where $\alpha_1$ is the curve defined in paragraph \ref{g}. On the other hand, from (\ref{simphi3}) and (\ref{acciong}) we have 
$(S_0^+)^{\ast}(g \Phi_3)= \overline{\frac{1}{g} \Phi_3}$.  Hence, taking into account (\ref{sinponer}) and that $\gamma_1=\alpha_1-(S_0^+)_{\ast}(\alpha_1)$ we get
$$ \int _{\gamma_1} g \Phi_3 = \int _{\alpha_1} g \Phi_3 -\overline{\int _{\alpha_1} \frac{1}{g} \Phi_3 } =2\: \ri \: {\rm Im} \left( \int _{\alpha_1} g \Phi_3 \right)\; .$$
Therefore, we have
\begin{equation} \label{dgamma1} d(a,\rho,\beta)= \frac{\ri}{2 \:\lambda}\int _{\gamma_1} g \Phi_3 = \frac{\ri}{2 \:\lambda}\int _{\widehat{\gamma}_1} g \Phi_3 \; ,\end{equation}
where $\widehat{\gamma}_1$ is a curve in $\overline{\cal N}$ homotopic to $\gamma_1$ that does not contain the point $1_+$. We can observe that when $\rho$ approaches $0$ the curve $\widehat{\gamma}_1$ is the boundary of a topological disk around the point $1_+$ and thereby we obtain
$$d(a,0,\beta)=  \frac{\ri}{2 \:\lambda} \lim_{\rho \to 0} \int _{\widehat{\gamma}_1} g \: \Phi_3 =
 \frac{-\pi}{\lambda}\:{\rm Residue}( \lim_{\rho \to 0} g \: \Phi_3,1_+) \; ,$$
Taking part \emph{a)} of Lemma \ref{lemaa} into account, we have that, when $\rho$ approaches $0$, $\frac{dg}{g}$ is a holomorphic one-form in a neighborhood of $1_+$ and furthermore $\lim_{\rho \to 0} g(1_+)=1$. On the other hand, we have that $\lim_{\rho \to 0} \Phi_3$ has a simple pole at $1_+$ and so
$$d(a,0,\beta)= \frac{-\pi}{\lambda} \: {\rm Residue}( \lim_{\rho \to 0} \Phi_3,1_+)=
\pi \frac{a (1+b^2)}{b (1+a^2)}>0 \; .$$
\noindent{\sc Proof of claim D.2 .} Taking into account the assertion \emph{b)} in Lemma \ref{lemaa}, the Gauss map can be easily computed in the case $\rho=\pi$ and we obtain
$$g(z)=\frac{z+\ri b}{\ri b z+1}\left( \frac{z-\ri a}{-\ri a z+1} \right)^{\beta} \; ,$$
where $b$ is given by part \emph{b)} in Lemma \ref{lemab} and we choose the branch of $z^\beta$ satisfying $1^\beta=1$. From here we infer that if $a \in [0,1[$ then $g$ has no poles in the curve $\alpha_1$. On the other hand, the expression of the one-form $\Phi_3$ in the case $\rho=\pi$ is given by
$$ \frac{1}{\lambda}\Phi_3(z)=-\frac{a}{b}\frac{(z+\ri b)(\ri b z+1)}{(z-\ri a)(-\ri a z+1)}\frac{dz}{z^2+1} \; .$$
Clearly, for $a \in [0,1[$, this one-form has a simple pole at $z=\pm \ri$ that are the extremes of the curve $\alpha_1$. Let us compute the sign of $-{\rm Im}\left({\rm Residue}\left(\frac{1}{\lambda}g \: \Phi_3,\pm \ri \right)\right)$. Then, we get
$$ {\rm Residue}\left(\frac{1}{\lambda}g \: \Phi_3,\ri \right)=\frac{a(1+b)^2 (1-a)^{\beta-1}}{2 b(1+a)^{\beta+1}}\: \re^{\beta \frac{\pi}{2}}\; , {\rm Residue}\left(\frac{1}{\lambda}g \: \Phi_3,-\ri \right)=-\frac{a(1-b)^2 (1+a)^{\beta-1}}{2 b(1-a)^{\beta+1}}\: \re^{-\beta \frac{\pi}{2}}$$
Thus we obtain $-{\rm Im}\left({\rm Residue}\left(\frac{1}{\lambda}g \: \Phi_3,\pm \ri \right)\right)<0$.
Note that if $a=1$, from statement \emph{b)} of Lemma \ref{lemab} we have that $b=1$. Therefore, in the case $a=1$ we have
$$ \frac{1}{\lambda} g \: \Phi_3= \frac{\ri \: dz}{(-\ri z+1)^{\beta}(z-\ri)^{2-\beta}} $$
and in this case we have a pole of order $2-\beta$ at $z=\ri$ and as before we deduce that $d(1,\rho,\beta)$ diverges to $-\infty$ when $\rho \to \pi$.

\noindent{\sc Proof of claim D.3 .} As was indicated in Remark \ref{paco} our examples for parameters $(0,\rho,\beta)$ coincide with López-Martín examples for parameters $n=\frac{2}{\beta +1}$, $r=-\cos(\frac{\rho}{\beta +1})$ studied in \cite{l-m2} and \cite{l-m1}. Furthermore, our function $d$ corresponds with function $f$ studied in Lemma 3 of \cite{l-m1}. From the analysis of this function developed by López and the second author in that paper follows that there exists $\rho_1(\beta) \in \: ]0,\pi[$ that verifies the conditions of {\sc Claim D.3.}
\vspace{0.3cm}\\
\noindent{\sc Proof of claim D.4 .} Along this proof we assume $a=1$ and $\rho \in \: ]0,\rho_0(\beta)]$. Using (\ref{simphi3}), (\ref{acciong}), (\ref{eqsd}) and that $\alpha_1=-S_{\ast}(\delta)+\delta$ we obtain
$$d(a,\rho,\beta)=-\frac{1}{\lambda}{\rm Im}\left( \int_1^{\re^{\ri \frac{\rho}{2}}}\left( g+\frac{1}{g}\right)\Phi_3 \right)\;. $$
Hence we can write $d(a,\rho,\beta)=I_1-I_2$, where
$$ I_1= \int_1^{\re^{\ri \frac{\rho}{2}}}\left( |g|+\frac{1}{|g|}\right)\cos(\arg(g)){\rm Im}\left(-\frac{1}{\lambda}\Phi_3 \right) = \int_0^\rho \left( |g(t)|+\frac{1}{|g(t)|}\right)\cos(\theta(t)){\rm Im}\left(-\frac{1}{\lambda}\Phi_3(\re^{\ri \frac{t}{2}}) \right) \, ,$$
$$ I_2= \int_1^{\re^{\ri \frac{\rho}{2}}}\left( |g|-\frac{1}{|g|}\right)\sin(\arg(g)){\rm Re}\left(\frac{1}{\lambda}\Phi_3 \right) = \int_0^\rho \left( |g(t)|-\frac{1}{|g(t)|}\right)\sin(\theta(t)){\rm Re}\left(\frac{1}{\lambda}\Phi_3(\re^{\ri \frac{t}{2}}) \right) \, ,$$ 
where $g(t)=g(\re^{\ri \frac{t}{2}})$ and $\theta(t)=\arg \left( g(\re^{\ri \frac{t}{2}})\right) \in [-\frac{\pi}{2},\frac{\pi}{2}]$. First of all we get
\begin{equation}\label{real} {\rm Re} \left( \tfrac{1}{\lambda}\Phi_3(\re^{\ri \frac{t}{2}}) \right)= \frac{2 \cos(\frac{\rho}{2})-c_2}{8 \cos(\frac{t}{2})^2}>0   \;. \end{equation}
Furthermore, since $h(1,\rho,\beta)<0$ for $\rho \in [0,\rho_0(\beta)[$ we deduce from (\ref{hdef}) that $2 \cos(\rho) - 2 \cos(\frac{\rho}{2})c_2>2$ and therefore we obtain
\begin{equation} \label{imag} {\rm Im}\left(-\tfrac{1}{\lambda}\Phi_3(\re^{\ri \frac{t}{2}}) \right)=\frac{-2 
\cos(t)+ 2 \cos(\rho)- 2 \cos(\frac{\rho}{2})c_2}{8 \cos(\frac{t}{2})^2\sqrt{2(\cos(t)-\cos(\rho))}}> \frac{
\sin(\frac{t}{2})^2}{2 \cos(\frac{t}{2})^2\sqrt{2(\cos(t)-\cos(\rho))}} \geq 0 \;. 
\end{equation}

Next, we will compute the expressions of $\theta(t)$ and $g(t)$. Regarding $\theta(t)$ we obtain 
\begin{equation} \label{theta} 0 \leq \theta(t)=\arctan \left( \frac{1-b^2}{1+b^2} \tan (\frac{t}{2})\right) \leq \frac{t}{2}\leq \frac{\pi}{2} \; . \end{equation}
On the other hand, $|g(t)|$ is given by
\begin{equation} \label{gt} |g(t)|=\exp \left( \int_1^{\re^{\ri \frac{t}{2}}} {\rm Re}\left( \frac{dg}{g} \right) \right)=\exp \left(  G(t,\rho,\beta) \right) \, ,\end{equation}
where $$G(t,\rho,\beta)=-\frac{1}{2  a_2}\int_0^{t} \left( \frac{\sqrt{2(\cos(s)- \cos(\rho))}}{b^2+\frac{1}{b^2}+2\cos(s)}-A_0 \frac{1}{\sqrt{2(\cos(s)- \cos(\rho))}} \right)ds \; .$$
For this function we obtain the following facts:
\begin{lemma}\label{lemaGt} The function $G$ above described satisfies:
\begin{enumerate}[a)]
\item $G(t,\rho,\beta) \geq 0 \, ,$
\item $G(t,\rho,\beta) \leq \frac{1-b^2}{\sqrt{b^4+1+2 b^2 \cos(\rho)}} G(\rho,\beta)$, where $G(\rho,\beta)=\log \left( \sqrt{\frac{b^2+1+2 b \sin(\frac{\rho}{2})}{b^2+1-2 b \sin(\frac{\rho}{2})} }\right) \; .$
\end{enumerate}
\end{lemma}
\begin{proof} First of all we note that $G(0,\rho,\beta)=G(\rho,\rho,\beta)=0$. Evidently, we have
$$ \frac{\partial G}{\partial t}(t,\rho,\beta)=-\frac{1}{2  a_2} \left( \frac{\sqrt{2(\cos(t)- \cos(\rho))}}{b^2+\frac{1}{b^2}+2\cos(t)}-A_0 \frac{1}{\sqrt{2(\cos(t)- \cos(\rho))}} \right) \,.$$
Thereby, if $t_0$ is a critical point of $G$ as a function of $t$ we have at this point
\begin{equation} 2\:(\cos(t_0)-\cos(\rho))-A_0 \left( 2 \cos(t_0)+b^2+\frac{1}{b^2} \right)=0\; . \label{critico}		          \end{equation}
Furthermore, the second derivative of $G$ respect to $t$ is
\begin{eqnarray*} \frac{\partial^2 G}{\partial t^2}(t,\rho,\beta) & = &- \frac{\sin(t)}{4 \sqrt{2} a_2 (b^2+\frac{1}{b^2}+2\cos(t))^2 (\cos(t)- \cos(\rho))^{\frac{3}{2}}}\left( -A_0(b^2+\frac{1}{b^2}+2\cos(t))^2 \right. \\  & - & \left. 2(\cos(t)- \cos(\rho))(b^2+\frac{1}{b^2}+2\cos(t)) +8(\cos(t)- \cos(\rho))^2
  \right)\end{eqnarray*}
Therefore, at the point $t_0$ we get
$$  \frac{\partial^2 G}{\partial t^2}(t_0,\rho,\beta)= - \frac{\sin(t_0)}{2 \sqrt{2} a_2  (\cos(t_0)- \cos(\rho))^{\frac{3}{2}}}A_0 (A_0-1) \; ,  $$
and since part \emph{c)} of Lemma \ref{lemaa} guarantees that $0 \leq A_0 \leq 1$, we deduce that there exists at most one critical point of $G$ as a function of $t$ and this must be a maximum. From here follows statement \emph{a)} of the lemma.

Now, we recall that $A_0 \geq 0$ and $\sqrt{2(\cos(s)- \cos(\rho))}\leq 2 \cos (\frac{s}{2})$. Then we have
$$G(t,\rho,\beta)\leq -\frac{1}{a_2}\int_0^{t} \frac{\cos(\frac{s}{2})}{b^2+\frac{1}{b^2}+2\cos(s)}ds=\frac{1-b^2}{\sqrt{b^4+1+2 b^2 \cos(\rho)}}\log \left(\sqrt{\frac{b^2+1+2 b \sin(\frac{t}{2})}{b^2+1-2 b \sin(\frac{t}{2})} } \right)\; .$$
To conclude the proof of \emph{b)} it is suffices to note that $\log \left(\sqrt{\frac{b^2+1+2 b \sin(\frac{t}{2})}{b^2+1-2 b \sin(\frac{t}{2})} } \right)$ is increasing as a function of $t$.
\end{proof}
Now we will return to the study of $d$. From (\ref{imag}) and (\ref{theta}) we get the following inequality for the integral $I_1 $
\begin{equation} \label{i1} I_1>\int_0^\rho \frac{\sin(\frac{t}{2})^2}{\cos(\frac{t}{2})\sqrt{2(\cos(t)- \cos(\rho))} }=\frac{\pi}{2} \left( \sec(\frac{\rho}{2})-1\right) \; ,\end{equation}
On the other hand, taking into account (\ref{real}), assertions \emph{a)} and \emph{b)} of Lemma \ref{lemaGt}, and that $\exp(x)-\exp(-x)$ and $\sin(x)$ are increasing functions we obtain
\begin{eqnarray*}&   I_2 & \leq \left(\exp(\tfrac{1-b^2}{\sqrt{b^4+1+2 b^2 \cos(\rho)}}G(\rho,\beta))-\exp(-\tfrac{1-b^2}{\sqrt{b^4+1+2 b^2 \cos(\rho)}}G(\rho,\beta))\right)\frac{2 \cos(\frac{\rho}{2})-c_2}{8} \int_0^t\frac{ \sin(\frac{s}{2})}{\cos(\frac{s}{2})^2}ds \\ &=&
\sinh \left( \tfrac{1-b^2}{\sqrt{b^4+1+2 b^2 \cos(\rho)}}G(\rho,\beta) \right)\frac{(2 \cos(\frac{\rho}{2})-c_2)}{2} \left( \sec(\frac{t}{2})-1 \right)\,. \end{eqnarray*}
As $\sec(x)$ is also increasing it follows
\begin{equation} \label{i2} I_2 \leq \sinh \left( \tfrac{1-b^2}{\sqrt{b^4+1+2 b^2 \cos(\rho)}}G(\rho,\beta) \right)\frac{(2 \cos(\frac{\rho}{2})-c_2)}{2} \left( \sec(\frac{\rho}{2})-1 \right)\,. \end{equation}
As a consequence of (\ref{i1}) and (\ref{i2}), we have
$$ d(1,\rho,\beta)= I_1-I_2 > \frac{1}{2}\left( \sec(\frac{\rho}{2})-1 \right)\left(\pi-(2 \cos(\frac{\rho}{2})-c_2)\sinh \left( \tfrac{1-b^2}{\sqrt{b^4+1+2 b^2 \cos(\rho)}}G(\rho,\beta) \right)\right) \,.$$
We denote by $H(\rho,\beta)=(2 \cos(\frac{\rho}{2})-c_2)\sinh \left( \tfrac{1-b^2}{\sqrt{b^4+1+2 b^2 \cos(\rho)}}G(\rho,\beta) \right)$. Our next objective is to prove that $\pi-H(\rho,\beta)$ is non negative.
We will distinguish several cases.

 Suppose first that $\rho \in [0,2\arcsin (\frac{\pi}{4})]$. Taking into account $\frac{1-b^2}{\sqrt{b^4+1+2 b^2 \cos(\rho)}}\leq 1$ and that $\sinh(x)$ is an increasing function we have
\begin{equation} \pi-H(\rho,\beta)\geq \pi-\frac{2\: b \sin(\frac{\rho}{2})(2 \cos(\frac{\rho}{2})-c_2)}{\sqrt{b^4+1+2 b^2 \cos(\rho)}}=\pi - \frac{4\: b \sin(\frac{\rho}{2}) \cos(\frac{\rho}{2})}{\sqrt{b^4+1+2 b^2 \cos(\rho)}}-2 \sin(\frac{\rho}{2}) \;. \label{final}\end{equation}
Consider now the function $k_1(\rho,\beta)=\frac{ b }{\sqrt{b^4+1+2 b^2 \cos(\rho)}}$. By deriving here we get
$$\frac{\partial k_1}{\partial \beta}(\rho,\beta)= \frac{1-b^4}{(b^4+1+2 b^2 \cos(\rho))^{\frac{3}{2}}} b_\beta \; .$$
From (\ref{bb}) we obtain that the preceding derivative is non negative. Therefore, taking part \emph{a)} of Lemma \ref{lemab} into account we get $k_1(\rho,\beta)\leq k_1(\rho,1) =\frac{1}{2 \cos (\frac{\rho}{2})}$. Thus, we can conclude
$$\pi-H(\rho,\beta)\geq \pi-4  \sin(\frac{\rho}{2}) \geq 0 \; .$$
Using Lemma \ref{lemaro} we know that if $ \beta \in \:[\beta_1,1]$, where $\beta_1= \frac{\pi}{2\arcsin (\frac{\pi}{4})}-1$, we are in the preceding case.  Thereby, we only have to study the cases where $\beta \in [0,\beta_1]$ and $\rho \in [2\arcsin (\frac{\pi}{4}),\pi]$. Observe that in the remainder cases $\cos(\rho) <0$.

Now, we assume $\beta \in [0,0.56]$ and $\rho \in [2\arcsin (\frac{\pi}{4}),\pi]$. Consider now the function $k_2(\rho,\beta)=\frac{ 1 }{\sqrt{b^4+1+2 b^2 \cos(\rho)}}$. By deriving here we get
$$\frac{\partial k_2}{\partial \beta}(\rho,\beta)= \frac{-2 b (b^2+\cos(\rho))}{(b^4+1+2 b^2 \cos(\rho))^{\frac{3}{2}}} \: b_\beta \; .$$
Recall that $b_\beta >0$, there is a critical point when $b^2=-\cos(\rho)$. In order to see the character of the critical point we compute the 
second derivative at a critical point and obtain 
$$\frac{\partial^2 k_2}{\partial \beta^2}(\rho,\beta)= \frac{4 \cos(\rho)}{\sin(\rho)^3} \: b_{\beta}<0 \; .$$
Hence, the critical point is a maximum and $k_2(\rho,\beta)\leq \frac{1}{\sin(\rho)}$. Then, from (\ref{final}) we have
\begin{equation} \label{intermedia}  \pi-H(\rho,\beta) \geq \pi -2\: b - 2 \sin(\frac{\rho}{2}) \; . \end{equation}
Using statement \emph{d)} of Lemma \ref{lemab} and Lemma \ref{lemaro} we have $b(1,\rho,\beta)\leq b(1,\frac{\pi}{\beta+1},\beta)$. If we compute the derivative of this function respect to $\beta$  and use (\ref{br}) and (\ref{bb}) we obtain
\begin{eqnarray*}\frac{\partial b}{\partial \beta}(1,\tfrac{\pi}{\beta+1},\beta) &=& -\tfrac{\pi}{(\beta+1)^2} b_\rho+b_\beta=\sqrt{b^4+1+2 b^2 \cos(\tfrac{\pi}{\beta+1})} \left( \frac{\pi}{(\beta+1)^2} \int_0^b \frac{t^2 \sin (\frac{\pi}{\beta+1})}{(t^4+1+2 t^2 \cos(\frac{\pi}{\beta+1})^{\frac{3}{2}}}dt \right. \\
& &  \left. +  \int_0^1 \frac{(1-t^2)^2+2t^2 \cos(\frac{\pi}{2(\beta+1)})^2(1-\frac{\pi \beta}{(\beta+1)^2}\tan(\frac{\pi}{2(\beta+1)})} {(t^4+1+2 t^2 \cos(\frac{\pi}{\beta+1})^{\frac{3}{2}}}dt \right)  \end{eqnarray*}
Let us prove $\frac{\partial b}{\partial \beta}(1,\tfrac{\pi}{\beta+1},\beta)\geq 0$. To see this it suffices to prove that $k(\beta)=\frac{\pi \beta}{(\beta+1)^2}\tan(\frac{\pi}{2(\beta+1)})\leq 1$. But this is very easy to see since $k$ is an increasing function and $k(1)=\frac{\pi}{8}<1$. 

Summarizing we have that $b(1,\frac{\pi}{\beta+1},\beta)\leq b(1,\frac{\pi}{1.56},0.56)$. Then substituting in (\ref{intermedia}) we get
$$ \pi-H(\rho,\beta) \geq \pi -2\: b(1,\frac{\pi}{1.56},0.56) - 2 >0$$

Suppose now that $\beta \in [0.56,\beta_1]$ and $\rho \in [2\arcsin (\frac{\pi}{4}),\pi]$. From this assumptions and Lemma \ref{lemaro}, we can assert that $\rho \in  [2\arcsin (\frac{\pi}{4}),\frac{\pi}{1.56}]$. Let us consider the function $k_3(\rho,\beta)=2 \cos(\frac{\rho}{2})+c_2=2 \cos(\frac{\rho}{2})+\frac{\sqrt{b^4+1+2 b^2 \cos(\rho)}}{b}$. It is not hard to see that $\frac{\partial k_3}{\partial \rho}$ and $\frac{\partial k_3}{\partial \beta}$ are non negative. Since $b_\beta$ and $b_\rho$ where non negative also (see (\ref{bb}) and part \emph{d)} in Lemma \ref{lemab}) we obtain that $k_3$ are a decreasing function on $\rho$ and $\beta$. Therefore, $k_3(\rho, \beta) \leq k_3(2\arcsin (\frac{\pi}{4}),0.56)< 3.04$. 

On the other hand, reasoning as before it is easy to see that $1-b^2$ is a decreasing function on $\rho$ and $\beta$. 
Then, taking into account the results obtained before for the function $k_2$ and  that $\rho \in  [2\arcsin (\frac{\pi}{4}),\frac{\pi}{1.56}]$ we get
$$\frac{ 1-b^2 }{\sqrt{b^4+1+2 b^2 \cos(\rho)}} \leq \frac{1-b(1,2\arcsin (\frac{\pi}{4}),0.56)^2  }{\sin(\frac{\pi}{1.56})}< 0.79 \; .$$
Finally, the same argument used with the preceding functions prove that $G(\rho,\beta)$ is an increasing function on $\rho$ and $\beta$ and so $G(\rho,\beta)\leq G(\frac{\pi}{1.56},\beta_1)<1.12$. Thus we have
$$ \pi-H(\rho,\beta) \geq \pi -3.04 \sinh(0.79 \cdot 1.12)=0.08>0 \;.$$

\noindent Departamento de Geometr\'{\i}a y Topolog\'{\i}a, Universidad de Granada, 18071 Granada, Spain \\ e-mails: {\tt lferrer@ugr.es, fmartin@ugr.es}


\begin{thebibliography}{99} \small

\bibitem{complex} R.B. Burckel. {\em An introduction to classical analysis. Vol. 1} (Birkhäuser, Basel, 1979).

\bibitem{CKMR} P. Collin, R. Kusner, W.H. Meeks, H. Rosenberg. {\em The topology, geometry and conformal structures of properly embedded minimal surfaces}, preprint.

\bibitem{farkas} H.M. Farkas, I. Kra. {\em Riemann Surfaces}. Springer-Verlag New York, 1992.

\bibitem{f-mp} L. Ferrer, F. Martín. {\em Properly embedded minimal disks bounded by non-compact polygonal lines.} Pacific J. Math. {\bf 214}(1), 55-88 (2004).

\bibitem{hkw} D. Hoffman, H. Karcher, F. Wei. {\em The singly-periodic genus-one helicoid.} Comment. Math. Hel.{\bf 74}, 248-279 (1999).

\bibitem{hww} D. Hoffman, M. Weber, M. Wolf. {\em An embedded genus-one helicoid.} Preprint (2004).

\bibitem{hw} D. Hoffman,  F. Wei. {\em Deforming the singly periodic genus-one helicoid.} Experiment. Math. {\bf 11} (2002), no. 2, 207--218.

\bibitem{jenkins-serrin} H. Jenkins, J. Serrin. {\em Variational problems of minimal
surface type.  II. Boundary value problems for the minimal surface equation.} Arch. Rat.
Mech. Anal., {\bf 21}, 321-342 (1966).

\bibitem{karcher} H. Karcher. {\em Construction of minimal surfaces.} Surveys in 
Geometry 1989/90, University of Tokyo 1989. Also: Vorlesungsreihe Nr. 1, SFB 256, Bonn, 1989.

\bibitem{l-m2} F. J. López, F. Martín. {\em A uniqueness theorem for properly embedded minimal surfaces bounded by straight lines.} J. Austr. Math. Soc. (Series A) {\bf 69}, 362-402 (2000).

\bibitem{l-m1} F. J. López, F. Martín. {\em Minimal surfaces in a wedge of a slab.} Comm. Anal. Geom. {\bf 9} no. 4,  683-723 (2001).

\bibitem{l-r-w} F. J. López, M. Ritoré, F. Wei. {\em A characterization of Riemann's minimal
surfaces.} J. Differential Geom., {\bf 47} No 2,  376-397 (1997).

\bibitem{m-r}  W.H. Meeks, H. Rosenberg. {\em The geometry and conformal structure of properly embedded minimal surfaces of finite topology in $\r^3$.} Invent. Math. {\bf 114}, 625-639 (1993).

\bibitem{weber} M. Weber. {\em The genus one helicoid is embedded}. 1999. Habilitationsschrift, Bonn.
\end{thebibliography}
\end{document}